\documentstyle[11pt]{article}

 \newcounter{abceqn}

 \newcounter{abcfig}

\newcommand{\HH}{{\mathcal H}}
\newcommand{\LL}{{\mathcal L}}
\newcommand{\q}{\vec{q}}
\newcommand{\B}{{\cal B}}
\newcommand{\E}{{\mathcal E}}
\newcommand{\V}{{\mathcal V}}

\newcommand{\vc}{\vec{c}}

\newcommand{\vth}{\vartheta}
\newcommand{\be}{\beta}

\newcommand{\om}{\omega}
\newcommand{\Ga}{\Gamma}
\newcommand{\Om}{\Omega}

\newcommand{\pa}{\partial}
\newcommand{\tF}{\tilde{F}}

\newcommand{\tq}{\tilde{q}}

\newcommand{\F}{{\cal{F}}}
\newcommand{\cS}{{\cal{S}}}

\newcommand{\vq}{\vec{q}}
\newcommand{\vQ}{\vec{Q}}

\newcommand{\e}{\varepsilon}
\newcommand{\U}{{\cal U}}

\newcommand{\A}{{\cal A}}
\newcommand{\k}{\kappa}
\newcommand{\ga}{\gamma}

\newcommand{\dl}{\delta}
\newcommand{\Dl}{\Delta}
\newcommand{\tDl}{\tilde{\Delta}}
\newcommand{\th}{\theta}
\newcommand{\ra}{\rightarrow}
\newcommand{\al}{\alpha}
\newcommand{\sg}{\sigma}
\newcommand{\Sg}{\Sigma}

\newcommand{\z}{\zeta}

\newcommand{\La}{\Lambda}
\newcommand{\la}{\lambda}
\newcommand{\bq}{\bar{q}}

\newcommand{\nid}{\noindent}

\newcommand{\W}{{\cal W}}
\newcommand{\lag}{\langle}
\newcommand{\rag}{\rangle}
\newcommand{\D}{{\cal D}}
\renewcommand{\theequation}{\thesection.\arabic{equation}}
\newcommand{\eqnsection}[1]{
	\section{#1}
	\setcounter{equation}{0}
	\renewcommand{\theequation}{\thesection.\arabic{equation}}
	\setcounter{figure}{0}
	\renewcommand{\thefigure}{\thesection.\arabic{figure}}
	\setcounter{remark}{0}
	\renewcommand{\theremark}{\thesection.\arabic{remark}}
	\setcounter{theorem}{0}
	\renewcommand{\thetheorem}{\thesection.\arabic{theorem}}
	\setcounter{lemma}{0}
	\renewcommand{\thelemma}{\thesection.\arabic{lemma}}
}

\title{\bf Homoclinic Tubes in Discrete Nonlinear Schr\"{o}dinger 
Equation Under Hamiltonian Perturbations}

\author{ \\ \\ \\ \\ Yanguang  (Charles) Li \\ \\  \\ \\ 
Department of Mathematics \\ \\ 
University of Missouri \\ \\ Columbia, MO 65211}

\date{\today}

\renewcommand{\theequation}{\thesection.\arabic{equation}}

\begin{document}
\bibliographystyle{plain}
\maketitle
\newpage
\begin{abstract}    
In this paper, we study the discrete cubic nonlinear Schr\"odinger lattice 
under Hamiltonian perturbations. First we develop a complete isospectral 
theory relevant to the hyperbolic structures of the lattice without 
perturbations. In particular, B\"{a}cklund-Darboux transformations are 
utilized to generate heteroclinic orbits and Melnikov vectors. Then we give 
coordinate-expressions for persistent invariant manifolds and Fenichel 
fibers for the perturbed lattice. Finally based upon the above 
machinery, existence of codimension 2 
transversal homoclinic tubes is established through a Melnikov type 
calculation and an implicit function argument. We also discuss symbolic 
dynamics of invariant tubes each of which consists of a doubly infinite 
sequence of curve segments when the lattice is four dimensional. 
Structures inside the asymptotic manifolds of the transversal homoclinic 
tubes are studied, special orbits, in particular homoclinic orbits and 
heteroclinic orbits when the lattice is four dimensional, are studied.

\vspace{0.5in}

\end{abstract}

\newtheorem{lemma}{Lemma}
\newtheorem{theorem}{Theorem}
\newtheorem{corollary}{Corollary}
\newtheorem{remark}{Remark}
\newtheorem{definition}{Definition}
\newtheorem{proposition}{Proposition}
\newtheorem{assumption}{Assumption}

\newpage
\tableofcontents



\newpage
\eqnsection{Introduction}

The concept of a homoclinic tube was introduced by Silnikov 
\cite{Sil68} in a study on the structure of the neighborhood 
of a homoclinic tube asymptotic to an invariant torus $\sg$ 
under a diffeomorphism $F$ in a finite dimensional phase space. 
The  asymptotic torus is of saddle type. The homoclinic tube 
consists of a doubly infinite sequence of tori $\{ \sg_j,\ j=0, 
\pm 1, \pm 2, \cdot \cdot \cdot \}$ in the transversal 
intersection of the stable and unstable manifolds of $\sg$, 
such that $\sg_{j+1} = F \circ \sg_j$ for any $j$. It is a 
generalization of the concept of a transversal homoclinic 
orbit when the points are replaced by tori. Silnikov obtained 
a similar theorem on the symbolic dynamics structures in the 
neighborhood of the homoclinic tube as Smale's theorem for a 
transversal homoclinic orbit \cite{Sil67} \cite{Pal88}.

We are interested in homoclinic tubes for several reasons: 1. 
Especially in high dimensions, dynamics inside each invariant 
tubes in the neighborhoods of homoclinic tubes are often 
chaotic too. We call such chaotic dynamics ``{\em{chaos in 
the small}}'', and the symbolic dynamics of the invariant 
tubes ``{\em{chaos in the large}}''. Such cascade structures 
are more important than the structures in a neighborhood of a 
homoclinic orbit, when high or infinite dimensional dynamical 
systems are studied. 2. Symbolic dynamics structures in the 
neighborhoods of homoclinic tubes are more observable than 
in the neighborhoods of homoclinic orbits in numerical and 
physical experiments due to the large dimensionality and 
the robustness of the homoclinic tubes. 3. When studying high or infinite 
dimensional Hamiltonian system (for example, the discrete 
NLS or NLS equations under Hamiltonian perturbations), each 
invariant tube contains both KAM tori and stochastic layers 
(chaos in the small). Thus, not only dynamics inside each 
stochastic layer is chaotic, all these stochastic layers also 
move chaotically under Poincar\'e maps.

In \cite{Li99b}, we proved the existence of transversal 
homoclinic tubes which are asymptotic to locally invariant 
center manifolds for the cubic nonlinear Schr\"{o}dinger equaton 
under Hamiltonian perturbations. Due to the locally invariant nature 
of the center manifolds, the symbolic dynamics in the neighborhood 
of the homoclinic tube is very difficult to establish. The locally 
invariant nature is very difficult to control in infinite dimensions. 
Therefore, we need to seek low dimensionality and hopefully can 
establish the above mentioned symbolic dynamics for the low dimensional 
systems which is the first step toward establishing such symbolic 
dynamics for PDEs. Then naturally we want to study the finite 
difference discretization of the cubic nonlinear Schr\"{o}dinger 
equation under Hamiltonian perturbations (perturbed discrete NLS). 
The dimension of the perturbed discrete NLS is $n$. As shown in 
this paper, when $n=4$, the center manifold is actually invariant 
and Silnikov's theorem implies the symbolic dynamics of line segments.

Denote by $W^{(c)}$ a normally hyperbolic locally invariant 
center manifold, by $W^{(cu)}$ and $W^{(cs)}$ the locally invariant
center-unstable and center-stable manifolds such that $W^{(c)}\subset 
W^{(cu)} \cap W^{(cs)}$, and by $F^t$ the evolution operator of the 
system. We call $\HH$ a transversal homoclinic tube if $\HH \subset 
W^{(cu)} \cap W^{(cs)}$, the intersection between $W^{(cu)}$ and 
$W^{(cs)}$ is transversal at $\HH$, and $\HH$ has the same dimension 
with $W^{(c)}$. Let $\Sg$ be 
an appropriate Poincar\'e section, and $P$ is the Poincar\'e map 
induced by the flow $F^t$; then $\HH \cap \Sg$ is called a 
transversal homoclinic tube under the Poincar\'e map $P$.

In \cite{LM97}, the discrete cubic nonlinear Schr\"odinger equation 
under dissipative perturbations is studied. Existence of a 
symmetric pair of homoclinic orbits is established through a 
Melnikov calculation and a geometric argument. In \cite{LW97a}, 
symbolic dynamics in the neighborhood of the symmetric pair of 
homoclinic orbits is established. In contrast to these previous 
works, the current work is a study on Hamiltonian perturbations 
and homoclinic tubes.

The paper is organized as follows: Section 2 is on the formulation 
of the problem, section 3 is on the isospectral theory for the 
discrete NLS equation, section 4 is on the coordinatization for invariant 
submanifolds, section 5 is on the existence of transversal 
homoclinic tubes, section 6 is on the symbolic dynamics of segments 
for the case when the perturbed discrete NLS equation is 4-dimensional,
section 7 is on structures inside the 
asymptotic manifolds of the transversal homoclinic tubes, 
and section 8 is the conclusion.

\newpage
\eqnsection{Formulation of the Problem}

Consider the discretized cubic nonlinear Schr{\"{o}}dinger equation 
under Hamiltonian perturbations (perturbed discrete NLS),
\\
\begin{equation}
i \dot{q}_n = \rho_n \pa H/ \pa \bq_n \ , \label{HDNLS}
\end{equation}
\\
where 
\[
H = H_0 +\e H_1 \ ,
\]
\[
H_0={1 \over h^2}\sum^{N-1}_{n=0}
\bigg \{\bq_n(q_{n+1}+q_{n-1})-{2 \over h^2}(1+\om^2 h^2)\ln \rho_n
\bigg \} \ .
\]
$\e$ is the perturbation parameter, $\e H_1$ is the Hamiltonian 
perturbation, $i = \sqrt{-1}$ is the imaginary unit, $\om$ is a positive 
parameter, and $\rho_n = 1 +h^2|q_n|^2$. $q_n$ 
satisfies the periodic and even boundary conditions,
\begin{equation}
q_{n+N}=q_n,\ \ \ \ q_{-n}=q_n \ ,  \label{bouc}
\end{equation}
where $N$ is a positive integer $N \geq 3$. (\ref{HDNLS}) is a 
$2(M+1)$-dimensional system, where $M=N/2$ ($N$ even) and 
$M=(N-1)/2$ ($N$ odd).
\begin{remark}
When $\e = 0$, $\sum^{N-1}_{n=0}\bigg \{\bq_n(q_{n+1}+q_{n-1})\bigg \}$ 
itself is also a constant of motion. This invariant, together with $H_0$,   
implies that $\sum^{N-1}_{n=0}\ln \rho_n$ is a constant of motion too. 
Therefore, 
\begin{equation}
D^2\equiv \prod^{N-1}_{n=0}\rho_n
\label{constD}
\end{equation}
is a constant of motion.
\end{remark}
The phase space is defined as 
\begin{eqnarray}
   \cS &\equiv& \bigg\{\vq=\left( \begin{array}{c} q \\ r 
\end{array} 
\right)\bigg| \ r=-\bar{q}, \ q=(q_0,q_1,...,q_{N-1})^T, \nonumber \\
             & & q_{n+N}=q_n, \ \ q_{N-n}=q_n\bigg\} \ . \label{phsp}
\end{eqnarray}
In $\cS$ (viewed as a vector space over the real numbers), we 
define the inner product, for any two points $\vq^{\ (1)}$ and $\vq^{\ (2)}$,
as follows:
\[
\lag \vq^{\ (1)}, \vq^{\ (2)} \rag = \sum_{n=0}^{N-1} ( q_n^{(1)} 
\overline{q_n^{(2)}} + \overline{q_n^{(1)}} q_n^{(2)}) \ . 
\]
\nid
And the norm of $\vq$ is defined as $\|\vq \ \|^2 = \lag \vq,\vq \rag$. 
We also use the notation $\vq_n = (q_n, r_n)^T$, where $r_n = - \bq_n$.

We will study transversal homoclinic tubes in this phase space $\cS$. 
When $\e =0$,
the unperturbed Hamiltonian system is the integrable discrete cubic nonlinear 
Schr{\"{o}}dinger equation,
\[
i \dot{q}_n = {1 \over h^2}[q_{n+1}-2 q_n +q_{n-1}] + |q_n|^2(q_{n+1}+q_{n-1})
-2 \om^2 q_n \ .
\]
The Hamiltonian perturbation term $H_1$ can be very general. In this paper, 
for concreteness, we will study the simple example,
\[
H_1={1 \over h^2}\sum^{N-1}_{n=0}\bigg \{ \al_1 (q_n +\bq_n)
+\al_2 (q_n^2 +\bq^2_n) \bigg \} \ln \rho_n\ ,
\]
where $\al_1$ and $\al_2$ are real parameters. The corresponding equation 
is,
\begin{eqnarray}
i \dot{q}_n &=& {1 \over h^2}[q_{n+1}-2 q_n +q_{n-1}] 
+ |q_n|^2(q_{n+1}+q_{n-1}) \nonumber \\
&-& 2 \om^2 q_n + \e \bigg \{ \bigg [\al_1 (q_n +\bq_n)+\al_2 (q_n^2 
+\bq^2_n) \bigg ] q_n \nonumber \\
&+& [\al_1 + 2 \al_2 \bq_n ]{\rho_n \over h^2}\ln \rho_n \bigg \}\ .
\label{PDNLS}
\end{eqnarray}
We will show that there exist transversal homoclinic tubes for this 
system in the phase space $\cS$ through Melnikov analysis and an implicit 
function argument. 
\begin{remark}
In the literature \cite{CEMS96}, for integrable systems under Hamiltonian 
perturbations, an interesting fact is that the Melnikov function sometimes 
is identically zero. Consider the NLS equation under Hamiltonian 
perturbation \cite{CEMS96}
\begin{equation}
iq_t = - q_{xx} -2 |q|^2q - \e q_{xxxx}\ ,
\label{fpNLS}
\end{equation}
with the Hamiltonian
\[
H=H_0 +\e H_1 = \int_0^1 [|q_x|^2 - |q|^4 -\e |q_{xx}|^2] dx\ .
\]
In \cite{CEMS96}, the Melnikov function $M$ is built with $H_0$,
\[
M= \int_{-\infty}^{\infty} \{ H_0,H_1 \} (Q) dt
\]
where $Q$ is a heteroclinic orbit generated through B\"{a}cklund-Darboux 
transformations, which is asymptotic to a periodic orbit independent of 
$x$ in $|t| \ra \infty$ limit. The reason why $M$ is identically zero is 
as follows
\[
M=-{1 \over \e} \int_{-\infty}^{\infty} {d H \over dt}\bigg |_Q dt =  
-{1 \over \e}\bigg [ H(\lim_{t \ra \infty}Q)- H(\lim_{t \ra -\infty}Q) \bigg ]
=0\ .
\]
From the study \cite{LM94}, one doubts that $H_0$ is the proper 
invariant to build the Melnikov function, since the center-unstable 
manifold associated with the unstable mode is really a level set of 
the invariant $F_1$ (\ref{coF}). Nevertheless, direct calculation shows that
\begin{eqnarray*}
{\dl F_1 \over \dl q}(Q) &=&\mbox{ linear combination of}\ \bigg \{ 
{\dl H_0 \over \dl q}(Q), {\dl H_0 \over \dl \bq}(Q), {\dl I \over \dl q}(Q),
{\dl I \over \dl \bq}(Q)\bigg \} \\
& & \mbox{ with coefficients only dependent upon}\ t\ ,
\end{eqnarray*}
where $I=\int_0^1 |q|^2 dx$ is still invariant under the perturbed 
flow (\ref{fpNLS}). This relation and the above Melnikov function do 
not imply that the Melnikov function built with $F_1$ is identically 
zero, but offer a hint that it is probably identically zero too. From 
some trial calculations for this paper, if a perturbation term still keeps
\[
I=\sum_{n=0}^{N-1}[\bq_n (q_{n+1} + q_{n-1})]
\]
invariant, the corresponding term of the Melnikov function built from $\tF_1$ 
(\ref{coF}) is identically zero. Such perturbation terms are for example,
\[
H_1 = {1 \over h^2}\sum_{n=0}^{N-1}\tilde{\al}_2 |q_n|^2 \ln \rho_n\ .
\]
When the Melnikov function is identically zero, second order Melnikov 
function is needed to measure the splitting between center-unstable 
and center-stable manifolds \cite{Li01d}.
\end{remark}

\newpage
\eqnsection{Isospectral Theory for Discrete NLS Equation}

In this section, we will study the isospectral theory for the discrete 
nonlinear Schr\"odinger equation
\begin{equation}
i \dot{q}_n = {1 \over h^2}[q_{n+1}-2 q_n +q_{n-1}] + |q_n|^2(q_{n+1}+q_{n-1})
-2 \om^2 q_n \ , \label{DNLS}
\end{equation}
where $i = \sqrt{-1}$, $q_n$'s are complex variables, $n \in Z$; under 
periodic and even boundary conditions,
\[
q_{n+N}=q_n\ , \ \ \ \ q_{-n}=q_n\ . 
\]
Since it is brand new, we will present 
the isospectral theory for the discrete NLS equation in details.

The integrability of the NLS equation (\ref{DNLS}) 
is proven with the use of the discretized Lax pair \cite{AL76}:
\begin{eqnarray}
\varphi_{n+1}&=&L^{(z)}_n\varphi_n, \label{Lax1} \\
\dot{\varphi}_n&=&B^{(z)}_n\varphi_n, \label{Lax2}
\end{eqnarray}
\noindent
where
\begin{eqnarray*}
L^{(z)}_n&\equiv&\left( \begin{array}{cc} 
                                            z& ihq_n \cr
                                            ih\bq_n & 1/z \cr
                                       \end{array} \right),
\\
\\
B^{(z)}_n&\equiv&{i\over h^2}\left( \begin{array}{cc}
  1-z^2+2i\la h-h^2q_n\bq_{n-1}+\om^2 h^2 & -izhq_n+(1/z)ihq_{n-1}  \\
  -izh\bq_{n-1}+(1/z)ih\bq_n& 1/z^2-1+2i\la h+h^2\bq_nq_{n-1}-\om^2 h^2
                                 \end{array} \right),
\end{eqnarray*}
\noindent
and where $z\equiv \exp(i\la h)$. Compatibility of the over determined 
system (\ref{Lax1},\ref{Lax2}) gives the ``Lax representation'' 
\[
\dot{L}_n=B_{n+1}L_n-L_nB_n
\]
of the discrete NLS equation (\ref{DNLS}). Focusing attention upon 
the discrete spatial flow (\ref{Lax1}), we let $Y^{(1)},
Y^{(2)}$ be the fundamental solutions of the ODE (\ref{Lax1}), i.e. solutions
with the initial conditions:
$$
           Y_0^{(1)}=\left( \begin{array}{c}
                                    1 \\ 0
                             \end{array} \right), \ 
           Y_0^{(2)}=\left( \begin{array}{c}
                                    0 \\ 1
                             \end{array} \right).
$$
\noindent
The $Floquet \ discriminant$:
\begin{equation}
\Delta : C \times \cS \mapsto C
\end{equation}
\noindent
is defined by
\begin{equation}
\Delta(z;\vec{q}) \equiv tr\{ M(N;z;\vec{q})\} \ ,
\end{equation}
\noindent
where $\cS$ is the phase space defined in (\ref{phsp}),
and $M(n;z;\vec{q}) = \ \mbox{columns}\{Y_n^{(1)},Y_n^{(2)}\}$ is
the fundamental matrix of (\ref{Lax1}). 
\begin{remark}  
$\Delta(z;\vec{q})$ is a constant of motion for the integrable 
system (\ref{DNLS}) for any $z \in C$ \cite{MO95}. Since 
$\Delta(z;\vec{q})$ is a meromorphic function in $z$ of degree ($+N,-N$), 
the Floquet discriminant $\Delta(z;\vec{q})$ acts as a generating 
function for $(M+1)$ functionally independent constants of motion, and is 
the key to the complete integrability of the system (\ref{DNLS}), where 
$M = N/2$ ($N$ even), $M=(N-1)/2$ ($N$ odd).
\end{remark}
\nid
The Floquet theory here is not standard as can be seen from the 
Wronskian relation: 
\[
W_N(\psi^+,\psi^-)=D^2 \ W_0(\psi^+,\psi^-), 
\]
where $D$ is defined in (\ref{constD}),
\[
W_n(\psi^+,\psi^-)\equiv \psi_n^{(+,1)}\psi_n^{(-,2)}-\psi_n^{(+,2)}
\psi_n^{(-,1)},
\]
$\psi^+$ and $\psi^-$ are any two solutions to the 
linear system (\ref{Lax1}).
In fact, $W_{n+1}(\psi^+,\psi^-)=\rho_n W_n(\psi^+,\psi^-)$. Due to 
this nonstandardness, modifications of the 
usual definitions of spectral quantities \cite{Li92} are required.
The Floquet spectrum is defined as the closure of the complex $\la$ for 
which there exists a bounded eigenfunction to the ODE (\ref{Lax1}). In 
terms of the Floquet discriminant $\Dl$, this is given by 
\[
\sg(L)= \bigg \{ z \in C \ \bigg | \ -2D \leq \Delta(z;\vec{q}) \leq 2D
\bigg \}\ .
\]
Periodic and antiperiodic points $z^s$ are defined by
$$
             \Delta(z^s;\vec{q})=\pm 2D.
$$
\noindent
A critical point $z^c$ is defined by the condition
$$
            {d\Delta \over dz}\bigg |_{(z^c;\vec{q})} = 0.
$$
\noindent
A multiple point $z^{m}$ is a critical point which is also a periodic 
or antiperiodic point. The $algebraic \ \ multiplicity$ of 
$z^{m}$ is defined as the order of the zero of $\Delta(z)\pm 2D$. Usually
it is $2$, but it can exceed $2$; when it does equal $2$, we call the 
multiple point a $double \ point$, and denote it by $z^{d}$. The $geometric
\ multiplicity$ of $z^{m}$ is defined as the dimension of the periodic 
(or antiperiodic) eigenspace of (\ref{Lax1}) at $z^{m}$, and is 
either $1$ or $2$. 

The normalized Floquet discriminant $\tilde{\Dl}$ is defined as
\[
\tilde{\Dl} = \Dl / D\ .
\]
\begin{remark}[Continuum Limit]
In the continuum limit (i.e. $h\ra 0$), the Hamiltonian has a limit in 
the manner: $H/h \ra H_c$, where $H_c$ is the Hamiltonian for 
NLS PDE, $H_c=i\int^1_0 \{q_x\bar{q}_x+2\om^2|q|^2-|q|^4 \}dx$.
The Lax pair (\ref{Lax1};\ref{Lax2}) also tends to the corresponding Lax
pair for NLS PDE with spectral parameter $\la$ ($z=e^{i\la h}$) \cite{LM94}. 
If $Q\equiv \max_n\{|q_n|\}$ is finite, then $\rho_n\ra 1$ as $h\ra 0$.
Therefore, $D^2\equiv \{\prod^{N-1}_{n=0}\rho_n\}\ra 1$
as $h\ra 0$. The nonstandard Floquet theory for the spatial 
part of the Lax pair (\ref{Lax1}) becomes the standard Floquet theory in the 
continuum limit.
\end{remark}

\subsection{Examples of Floquet Spectra \label{EFS}}

Consider the uniform solution: $q_n = q_c$, $\forall n$ 
\[
q_c(t) = a \exp \bigg \{ -i [ 2(a^2 - \om^2)t - \ga] \bigg \}\ .
\]
The corresponding Bloch functions of the Lax pair are given by:
\begin{eqnarray*}
\psi^+_n &=& (\sqrt{\rho} e^{i\be})^n e^{\Om_+ t} \left ( 
\begin{array}{c} ({1 \over z} - \sqrt{\rho}e^{i\be}) \exp 
\{ -i [ (a^2 - \om^2)t - \ga/2]\} \cr \cr -iha  \exp \{ i [ (a^2 - \om^2)t 
- \ga/2]\} \cr  \end{array} \right ) \ , \\ 
\psi^-_n &=& (\sqrt{\rho} e^{-i\be})^n e^{\Om_- t} \left ( 
\begin{array}{c} -iha  \exp \{ -i [ (a^2 - \om^2)t 
- \ga/2]\} \cr \cr (z - \sqrt{\rho}e^{-i\be}) \exp 
\{ i [ (a^2 - \om^2)t - \ga/2]\} \cr \end{array} \right ) \ , 
\end{eqnarray*}
where
\[
z=\sqrt{\rho} \cos \be + \sqrt{\rho \cos^2 \be - 1}\ ,\ \ 
\rho = 1 + h^2 a^2\ ,
\]
\[
\Om_+ = {i \over h^2} \bigg \{ ({1 \over z} - z ) \sqrt{\rho} e^{i\be} 
+ i 2 \la h \bigg \} \ ,
\]
\[
\Om_- = {i \over h^2} \bigg \{ ({1 \over z} - z ) \sqrt{\rho} e^{-i\be} 
+ i 2 \la h \bigg \} \ .
\]
The Floquet discriminant is given by:
\[
\Dl = 2 D \cos (N\be)\ ,
\]
where $D= \rho^{N/2}$. Thus the Floquet spectra are given by:
\[
-1 \leq \cos (N\be) \leq 1 \ .
\]
See figure \ref{efl} for an illustration of the Floquet spectra.
\begin{figure}
\vspace{3.0in}
\caption{An illustration of the Floquet spectra for the uniform potential 
$q_c$ (The number and locations of the spectral points can vary).}
\label{efl}
\end{figure}
Periodic and antiperiodic points are given by:
\[
z^{(s)}_m=\sqrt{\rho} \cos {m \over N} \pi + 
\sqrt{\rho \cos^2 {m \over N} \pi - 1}\ ,
\]
where $z^{(s)}_m$ is a periodic point when $m$ is even, and $z^{(s)}_m$ is an 
antiperiodic point when $m$ is odd. See figure \ref{efl} for an 
illustration of the periodic and antiperiodic points.
The following facts are obvious,
\begin{itemize}
\item If $N$ is odd, all the periodic and antiperiodic points are on 
the real axis for sufficiently large $|q_c|$.
\item If $N$ is even, except the two points $z=\pm i$, all other 
periodic and antiperiodic points are on 
the real axis for sufficiently large $|q_c|$.
\end{itemize}
Derivatives of the Floquet discriminant $\Dl$ with respect to $z$ are 
given by:
\[
d\Dl / dz = 2ND \sin (N\be) [z \sqrt{\rho} \sin \be ]^{-1} \sqrt{
\rho \cos^2 \be -1} \ ,
\]
\begin{eqnarray*}
d^2\Dl / dz^2 &=& - 2ND [ \rho z^2 \sin^3 \be ]^{-1} \bigg [
N \cos(N\be) \sin \be \ (\rho \cos^2 \be - 1) \\
&+&(1 - \rho) \cos \be \sin (N\be) 
+ \sqrt{\rho}\sin (N\be)\sin^2 \be \sqrt{\rho \cos^2 \be - 1} \bigg ] \ .
\end{eqnarray*}
The critical points are given by:
\[
\be = { m \over N} \pi \ \ ( \be \neq 0, \pi), \ \ \mbox{or}\ \cos^2 
\be = 1/\rho \ .
\]
See figure \ref{efl} for an 
illustration of the critical points and double points. There can be multiple 
points with algebraic multiplicity greater than $2$. For example, when 
two symmetric double points on the circle collide at the intersection points 
$z = \pm 1$, we have multiple points of algebraic multiplicity $4$. In such 
case, $\rho \cos^2 \be - 1 = 0$ and $\sin (N\be) = 0$; then $d^2\Dl / dz^2 
= 0$, at $z = \pm 1$.

\subsection{An Important Sequence of Invariants}

\begin{definition}
The sequence of invariants $\tF_j$ is defined as:
\begin{equation}
\tF_j(\vq) = \tilde{\Dl}(z^{(c)}_j(\vq);\vq)\ .
\label{coF}
\end{equation}
\end{definition}
These invariants $\tF_j$'s are perfect candidate for building Melnikov 
functions. The Melnikov vectors are given by the gradients of these 
invariants.
\begin{lemma}
Let $z^{(c)}_j(\vq)$ be a simple critical point; then
\begin{equation}
{\dl \tF_j \over \dl \vq_n}(\vq) = {\dl \tilde{\Dl} \over \dl \vq_n}
(z^{(c)}_j(\vq);\vq)\ .
\label{derf}
\end{equation}
\begin{equation}
{\dl \tilde{\Dl} \over \dl \vq_n}(z;\vq) = {i h (\z -\z^{-1}) \over 
2 W_{n+1}} \left (\begin{array}{c} \psi^{(+,2)}_{n+1}\psi^{(-,2)}_{n}+ 
\psi^{(+,2)}_{n}\psi^{(-,2)}_{n+1} \cr \cr \psi^{(+,1)}_{n+1}\psi^{(-,1)}_{n}+ 
\psi^{(+,1)}_{n}\psi^{(-,1)}_{n+1} \cr \end{array} \right )\ , 
\label{fidr}
\end{equation}
where $\psi_n^\pm = (\psi_n^{(\pm,1)}, \psi_n^{(\pm,2)})^T$ are two 
Bloch functions of the Lax pair (\ref{Lax1},\ref{Lax2}), such that 
\[
\psi^\pm_n = D^{n/N} \z^{\pm n/N} \tilde{\psi}^\pm_n\ ,
\]
where $\tilde{\psi}^\pm_n$ are periodic in $n$ with period $N$,
$W_n = \ \mbox{det}\ ( \psi^+_n,\psi^-_n )$.
\end{lemma}

Proof: By the definition of critical points,
\[
\tDl'(z^{(c)}_j(\vq);\vq) = 0\ .
\]
Differentiating this equation, we have 
\[
{\dl z^{(c)}_j \over \dl \vq_n} = - { 1 \over \tDl''}{\dl \tDl' \over \dl 
\vq_n} \ .
\]
Since $z^{(c)}_j(\vq)$ is a simple critical point, $z^{(c)}_j$ is a 
differentiable function. Thus
\begin{eqnarray*}
{\dl \tF_j \over \dl \vq_n} &=& {\dl \tilde{\Dl} \over \dl \vq_n}\bigg |_{
z=z^{(c)}_j} + {\pa \tDl \over \pa z}\bigg |_{z=z^{(c)}_j}{\dl z^{(c)}_j 
\over \dl \vq_n} \\ 
&=& {\dl \tilde{\Dl} \over \dl \vq_n}\bigg |_{z=z^{(c)}_j}\ .
\end{eqnarray*}
This proves equation (\ref{derf}). Next we derive the formula (\ref{fidr}).
Let $M_n$ be the fundamental matrix to the system (\ref{Lax1}), i.e. 
the matrix solution to (\ref{Lax1}) with initial condition $M_0$ being 
a 2 x 2 identity matrix. Variation of $\vq$ leads to the variational equation 
for the variation of $M_n$ at fixed $z$,
\begin{eqnarray*}
\dl M_{n+1} &=& \left ( \begin{array}{lr} z & ihq_n \cr \cr ih\bq_n & 1/z \cr 
\end{array}\right ) \dl M_n + \left ( \begin{array}{lr} 0 & ih\ \dl q_n \cr 
\cr ih\ \overline{\dl q_n} & 0 \cr 
\end{array}\right )  M_n \ , \\ \\ 
\dl M_0 &=& 0\ .
\end{eqnarray*}
Let $\dl M_n = M_n A_n$, where $A_n$ is a 2 x 2 matrix to be determined, we 
have
\begin{eqnarray}
A_{n+1} - A_n &=& M_{n+1}^{-1} \dl U_n M_n\ , \nonumber \\
\label{varA} \\
A_0 &=& 0 \ , \nonumber 
\end{eqnarray}
where
\[
\dl U_n = \left ( \begin{array}{lr} 0 & ih\ \dl q_n \cr \cr 
ih\ \overline{\dl q_n} & 0\cr \end{array}\right ) \ .
\]
Solving the system (\ref{varA}), we have
\begin{eqnarray*}
\dl M_{n} &=& M_n \bigg [ \sum^n_{j=1} M_j^{-1} \dl U_{j-1} 
M_{j-1} \bigg ]\ ,\\
\dl M_0 &=& 0\ .
\end{eqnarray*}
Then,
\[
\dl \Dl(z,\vq) = \ \mbox{trace}\ \bigg \{ M_N [ \sum^N_{j=1} M_j^{-1} 
\dl U_{j-1} M_{j-1}] \bigg \}\ .
\]
Thus,
\begin{eqnarray}
{\dl \Dl \over \dl q_n} &=& i h \ \mbox{trace}\ \bigg \{ M_{n+1}^{-1} 
\left ( \begin{array}{lr} 0 & 1 \cr 0&0 \cr \end{array} \right ) M_n M_N 
\bigg \}\ , \label{fir1} \\ 
{\dl \Dl \over \dl r_n} &=& -i h \ \mbox{trace}\ \bigg \{ M_{n+1}^{-1} 
\left ( \begin{array}{lr} 0 & 0 \cr 1&0 \cr \end{array} \right ) M_n M_N 
\bigg \}\ , \label{fir2}
\end{eqnarray}
where $r_n = - \bq_n$. Let $\psi^+$ and $\psi^-$ be two Bloch functions for 
the discrete Lax pairs (\ref{Lax1},\ref{Lax2}),
\begin{equation}
\psi^{\pm}_n = D^{n/N} \z^{\pm n/N} \tilde{\psi}^{\pm}_n\ ,
\label{blocfuc}
\end{equation}
where $\tilde{\psi}^{\pm}_n$ are periodic in $n$ with period $N$. Let $B_n$ 
be the 2 x 2 matrix with $\psi_n^+$ and $\psi_n^-$ as the column vectors,
\begin{equation}
B_n = ( \psi_n^+, \psi_n^- )\ . \label{blocm}
\end{equation}
Then
\begin{equation}
M_n = B_n B_0^{-1}\ , \ \ \ M_N = B_0 \left ( \begin{array}{lr} 
D\z & 0 \cr 0 & D\z^{-1} \cr \end{array} \right ) B_0^{-1}\ . 
\label{blocr}
\end{equation}
Substitute the representations (\ref{blocm}) and (\ref{blocr}) into 
(\ref{fir1},\ref{fir2}), we have
\[
{\dl \Dl \over \dl \vq_n} = \tDl {\dl D \over \dl \vq_n} +
{i D h (\z -\z^{-1}) \over 
2 W_{n+1}} \left (\begin{array}{c} \psi^{(+,2)}_{n+1}\psi^{(-,2)}_{n}+ 
\psi^{(+,2)}_{n}\psi^{(-,2)}_{n+1} \cr \cr \psi^{(+,1)}_{n+1}\psi^{(-,1)}_{n}+ 
\psi^{(+,1)}_{n}\psi^{(-,1)}_{n+1} \cr \end{array} \right )\ , 
\]
where
\[
W_n = \ \mbox{det}\ B_n\ , \ \ \ \psi^{\pm}_n = (\psi^{(\pm,1)}_n,
\psi^{(\pm,2)}_n)^T\ . 
\]
Thus,
\[
{\dl \tDl \over \dl \vq_n} = {i h (\z -\z^{-1}) \over 
2 W_{n+1}} \left (\begin{array}{c} \psi^{(+,2)}_{n+1}\psi^{(-,2)}_{n}+ 
\psi^{(+,2)}_{n}\psi^{(-,2)}_{n+1} \cr \cr \psi^{(+,1)}_{n+1}\psi^{(-,1)}_{n}+ 
\psi^{(+,1)}_{n}\psi^{(-,1)}_{n+1} \cr \end{array} \right )\ , 
\]
which is the formula (\ref{fidr}). The lemma is proved. $\Box$

\subsection{B\"acklund-Darboux Transformations}
The hyperbolic structure and homoclinic orbits for (\ref{DNLS}) are 
constructed through the B\"acklund-Darboux transformations, which were
built in \cite{Li92}. First, we present the B\"acklund-Darboux transformations.
Then, we show how to construct homoclinic orbits.

Fix a solution $q_n(t)$ of the 
system (\ref{DNLS}), for which the linear operator $L_n$ has a double 
point $z^d$ of geometric multiplicity 2, which is not on the unit circle.
We denote two linearly independent solutions (Bloch functions) of the 
discrete Lax pair (\ref{Lax1};\ref{Lax2}) at $z=z^d$ by $(\phi_n^+,\phi_n^-)$.
Thus, a general solution of the discrete Lax pair (\ref{Lax1};\ref{Lax2}) 
at $(q_n(t),z^d)$ is given by
\[
\phi_n(t; z^d, c^+, c^-)=c^+ \phi_n^+ + c^- \phi_n^-,
\]
\nid
where $c^+$ and $c^-$ are complex parameters. We use $\phi_n$ 
to define a transformation matrix $\Ga_n$ by
\[
\Ga_n=\left(\begin{array}{cc} z+(1/z)a_n & b_n \cr c_n &-1/z+z d_n \cr
            \end{array} \right),
\]
\nid
where,
\begin{eqnarray*}
a_n &=& {z^d \over (\bar{z}^d)^2\Dl_n}\bigg [|\phi_{n2}|^2+|z^d|^2|\phi_{n1}|^2
    \bigg ],\\
d_n &=& -{1 \over z^d\Dl_n}\bigg [|\phi_{n2}|^2+|z^d|^2|\phi_{n1}|^2
    \bigg ],\\
b_n &=& {|z^d|^4-1 \over (\bar{z}^d)^2\Dl_n}\phi_{n1}\bar{\phi}_{n2}, \\
c_n &=& {|z^d|^4-1 \over z^d\bar{z}^d\Dl_n}\bar{\phi}_{n1}\phi_{n2}, \\
\Dl_n &=& -{1 \over \bar{z}^d}\bigg [|\phi_{n1}|^2+|z^d|^2|\phi_{n2}|^2
    \bigg ].
\end{eqnarray*}
\nid
From these formulae, we see that
\[
\bar{a}_n=-d_n,\ \ \bar{b}_n=c_n.
\]
\nid
Then we define $Q_n$ and $\Psi_n$ by
\begin{equation}
Q_n\equiv {i\over h}b_{n+1}-a_{n+1}q_n
\label{BD1}
\end{equation}
\nid
and
\begin{equation}
\Psi_n(t;z)\equiv \Ga_n(z;z^d;\phi_n)\psi_n(t;z)
\label{BD2}
\end{equation}
\nid
where $\psi_n$ solves the discrete Lax pair (\ref{Lax1};\ref{Lax2}) 
at $(q_n(t),z)$. Formulas (\ref{BD1}) and (\ref{BD2}) are the 
B\"acklund-Darboux transformations for the potential and eigenfunctions, 
respectively. We have the following theorem \cite{Li92}.
\begin{theorem}[B\"acklund-Darboux Transformations]
Let $q_n(t)$ denote a solution of the system (\ref{DNLS}), for which 
the linear operator $L_n$ has a double point $z^d$ of geometric multiplicity 2, 
which is not on the unit circle and which is associated with an 
instability. We denote two linearly independent 
solutions (Bloch functions) of the discrete Lax pair (\ref{Lax1};\ref{Lax2}) 
at $(q_n, z^d)$ by $(\phi_n^+,\phi_n^-)$. We define $Q_n(t)$ and $\Psi_n(t;z)$ 
by (\ref{BD1}) and (\ref{BD2}). Then
\begin{enumerate}
\item $Q_n(t)$ is also a solution of the system (\ref{DNLS}). (The eveness
of $Q_n$ can be obtained by choosing the complex B\"acklund parameter 
$c^+/c^-$ to lie on a certain curve, as shown in the example below.)
\item $\Psi_n(t;z)$ solves the discrete Lax pair (\ref{Lax1};\ref{Lax2}) 
at $(Q_n(t),z)$.
\item $\Dl(z;Q_n)=\Dl(z;q_n)$, for all $z\in C$.
\item $Q_n(t)$ is homoclinic to $q_n(t)$ in the sense that $Q_n(t) \ra 
e^{i\th_{\pm}}\ q_n(t)$, exponentially as $\exp (-\sg |t|)$ as $t \ra 
\pm \infty$. Here $\th_{\pm}$ are the phase shifts, $\sg$ is a 
nonvanishing growth rate associated to the double point $z^d$, and 
explicit formulas can be developed for this growth rate and 
for the phase shifts $\th_{\pm}$.
\end{enumerate}
\label{Backlund}
\end{theorem}
This theorem is quite general, constructing homoclinic solutions from 
a wide class of starting solutions $q_n(t)$. Its proof is by direct 
verification \cite{Li92}. 

\subsection{Heteroclinic Orbits and Melnikov Vectors}

Through B\"acklund-Darboux transformations, heteroclinic orbits $Q_n$ 
are generated, with the formula expression (\ref{BD1})
\begin{eqnarray*}
Q_n &=& i h^{-1} b_{n+1} - a_{n+1} q_n \\
    &=& \bigg [ \overline{z^{(d)}} (|\phi^{(1)}_{n+1}|^2 + |z^{(d)}|^2 
|\phi^{(2)}_{n+1}|^2) \bigg ]^{-1} \times \\
& & \bigg [ ih^{-1}(1-|z^{(d)}|^4)\phi^{(1)}_{n+1}\overline{\phi^{(2)}_{n+1}} 
+z^{(d)} q_n (|\phi^{(2)}_{n+1}|^2 + |z^{(d)}|^2 
|\phi^{(1)}_{n+1}|^2) \bigg ]\ .
\end{eqnarray*}
The Melnikov vector field located on this heteroclinic orbit is given by
\begin{equation}
{\dl \tDl(z^{(d)};\vQ) \over \dl \vQ_n} = K {W_n \over E_n A_{n+1}} 
\left ( \begin{array}{c} [z^{(d)}]^{-2} \ \overline{\phi^{(1)}_n} 
\ \overline{\phi^{(1)}_{n+1}} \cr \cr [\overline{z^{(d)}}]^{-2} 
\ \overline{\phi^{(2)}_n} \ \overline{\phi^{(2)}_{n+1}} \cr \end{array} 
\right )\ , \label{fibd}
\end{equation}
where 
\[
\phi_n = (\phi^{(1)}_n, \phi^{(2)}_n)^T= c_+ \psi_n^+ + c_- \psi_n^- \ ,
\]
\[
W_n = \left | \begin{array}{lr} \psi_n^+ & \psi_n^-  \cr \end{array} 
\right |\ ,
\]
\[
E_n = |\phi^{(1)}_n|^2 + |z^{(d)}|^2|\phi^{(2)}_n|^2\ ,
\]
\[
A_n = |\phi^{(2)}_n|^2 + |z^{(d)}|^2|\phi^{(1)}_n|^2\ ,
\]
\[
K=-{ihc_+c_- \over 2}|z^{(d)}|^4 (|z^{(d)}|^4 - 1)[\overline{z^{(d)}}]^{-1} 
\sqrt{ \tDl(z^{(d)};\vq) \tDl''(z^{(d)};\vq)}\ ,
\]
where
\[
\tDl''(z^{(d)};\vq) = {\pa^2 \tDl(z^{(d)};\vq) \over \pa z^2}\ .
\]

\subsubsection{An Important Example}

In this subsection, we start with the uniform solution of (\ref{DNLS})
\begin{equation}
q_n=q_c \ ,\ \forall n; \ \ \ \ q_c=a\exp \bigg \{-i[2(a^2-\om^2)t - \ga]
\bigg \}\ . 
\label{ucsl}
\end{equation}
We choose the amplitude $a$ in the range
\begin{eqnarray}
& & N\tan{\pi \over N}< a <N\tan{2\pi \over N}\ ,\ \ \ \mbox{when}\ N>3 \ ,
\nonumber \\ \label{constr} \\
& & 3\tan{\pi \over 3}< a < \infty\ ,\ \ \ \mbox{when}\ N=3 \ ;\nonumber
\end{eqnarray}
so that there is only one set of quadruplets of double points 
which are not on the unit circle, and denote one of them by $z=z_1^{(d)}=
z_1^{(c)}$ which corresponds to $\be = \pi / N$ (see subsection \ref{EFS} and 
figure \ref{efl}). The heteroclinic orbit $Q_n$ is given by
\begin{equation}
Q_n = q_c (\hat{E}_{n+1})^{-1} \bigg [ \hat{A}_{n+1} - 2 \cos \be 
\sqrt{\rho \cos^2 \be -1}\hat{B}_{n+1} \bigg ]\ ,
\label{hetorb}
\end{equation}
and the Melnikov vectors evaluated on this heteroclinic orbit are given 
by
\begin{equation}
{\dl \tF_1 \over \dl \vQ_n} = \hat{K} \bigg [ \hat{E}_n \hat{A}_{n+1}
\bigg ]^{-1} \ \mbox{sech}\ [2\mu t +2p]\left ( \begin{array}{c} 
\hat{X}^{(1)}_n \cr - \hat{X}^{(2)}_n \cr \end{array} \right ) \ ,
\label{melv}
\end{equation}
where
\[
\hat{E}_n = ha\cos \be +\sqrt{\rho \cos^2 \be - 1} \ \mbox{sech}\ [
2\mu t +2p] \cos [(2n-1)\be +\vth ]\ ,
\]
\[
\hat{A}_{n+1} = ha\cos \be +\sqrt{\rho \cos^2 \be - 1} \ \mbox{sech}\ [
2\mu t +2p] \cos [(2n+3)\be +\vth ]\ ,
\]
\[
\hat{B}_{n+1} = \cos \varphi + i \sin \varphi \tanh [ 2\mu t +2p]
+ \ \mbox{sech}\ [ 2\mu t +2p] \cos [ 2(n+1)\be +\vth ]\ ,
\]
\begin{eqnarray*}
\hat{X}^{(1)}_n &=& \bigg [ \cos \be \ \mbox{sech}\ [ 2\mu t +2p] 
+\cos [(2n+1)\be +\vth +\varphi] \\ 
& & - i \tanh [ 2\mu t +2p] 
\sin [(2n+1)\be +\vth +\varphi] \bigg ] e^{i2\th(t)}\ ,
\end{eqnarray*}
\begin{eqnarray*}
\hat{X}^{(2)}_n &=& \bigg [ \cos \be \ \mbox{sech}\ [ 2\mu t +2p] 
+\cos [(2n+1)\be +\vth -\varphi] \\
& & - i \tanh [ 2\mu t +2p] 
\sin [(2n+1)\be +\vth -\varphi] \bigg ] e^{-i2\th(t)}\ ,
\end{eqnarray*}
\[
\hat{K} = -2Nh^2a(1-z^4) [8 \rho^{3/2}z^2]^{-1}\sqrt{\rho \cos^2 \be -1}\ ,
\]
\[
\be = \pi / N\ , \ \ \rho = 1+h^2 a^2\ , \ \ \mu = 2h^{-2} \sqrt{\rho} 
\sin \be \sqrt{\rho \cos^2 \be -1}\ ,
\]
\[
h=1/N,\ \ c_+/c_- = i e^{2p} e^{i\vth}\ , \ \ \vth \in [0,2\pi]\ , \ \ 
p \in (-\infty, \infty)\ ,
\]
\[
z=\sqrt{\rho}\cos \be +\sqrt{\rho \cos^2 \be -1}\ , \ \ 
\th(t)=(a^2-\om^2)t - \ga/2\ ,
\]
\[
\sqrt{\rho \cos^2 \be -1} + i \sqrt{\rho} \sin \be = ha e^{i \varphi}\ ,
\]
where $\varphi=\sin^{-1} [\sqrt{\rho}(ha)^{-1} \sin \be ],
\ \ \varphi \in (0, \pi/2)$.

Next we study the ``evenness'' condition: $Q_{-n} = Q_n$. It turns out 
that the choices $\vth = - \be\ , \ - \be +\pi$ in the formula of 
$Q_n$ lead to the evenness of $Q_n$ in $n$. In terms of figure eight 
structure of $Q_n$, $\vth = - \be$ corresponds to one ear of the 
figure eight, and $\vth = - \be +\pi$ corresponds to the other 
ear. The even formula for $Q_n$ is given by,
\begin{equation}
Q_n = q_c \bigg [ \Ga / \La_n -1 \bigg ]\ , \label{ehetorb}
\end{equation}
where
\[
\Ga = 1-\cos 2 \varphi - i \sin 2 \varphi \tanh [ 2 \mu t + 2p]\ ,
\]
\[
\La_n = 1 \pm \cos \varphi [\cos \be ]^{-1}\ \mbox{sech}[2 \mu t + 2p] 
\cos [2n\be]\ ,
\]
where (`+' corresponds to $\vth = -\be$).
The Melnikov vectors evaluated on these heteroclinic orbits are not 
necessarily even and are in fact not even in $n$. For the purpose of 
calculating the Melnikov functions, only the even parts of the 
Melnikov vectors are needed, which are given by
\begin{equation}
{ \dl \tF_1 \over \dl \vQ_n}\bigg |_{\mbox{even}} = \hat{K}^{(e)} 
\ \mbox{sech}[2\mu t + 2p] [\Pi_n]^{-1}\left ( \begin{array}{c} 
\hat{X}^{(1,e)}_n \cr - \hat{X}^{(2,e)}_n \cr \end{array} \right ) \ ,
\label{emelv}
\end{equation}
where
\[
\hat{K}^{(e)}= -2N (1-z^4) [8a\rho^{3/2} z^2]^{-1} 
\sqrt{\rho \cos^2\be - 1}\ ,
\]
\begin{eqnarray*}
\Pi_n &=& \bigg [ \cos \be \pm \cos \varphi \ \mbox{sech}[2\mu t +2p] 
\cos[2(n-1)\be]\bigg ] \times \\
& &\bigg [ \cos \be \pm \cos \varphi \ \mbox{sech}[2\mu t +2p] 
\cos[2(n+1)\be]\bigg ]\ ,
\end{eqnarray*}
\begin{eqnarray*}
\hat{X}^{(1,e)}_n &=& \bigg [ \cos \be \ \mbox{sech}\ [ 2\mu t +2p] 
\pm (\cos \varphi \\
& & -i \sin \varphi \tanh [ 2\mu t +2p]) \cos [2n\be]\bigg ] e^{i2\th(t)}\ ,
\end{eqnarray*} 
\begin{eqnarray*}
\hat{X}^{(2,e)}_n &=& \bigg [ \cos \be \ \mbox{sech}\ [ 2\mu t +2p] 
\pm (\cos \varphi \\
& & +i \sin \varphi \tanh [ 2\mu t +2p]) \cos [2n\be]\bigg ] e^{-i2\th(t)}\ .
\end{eqnarray*}

The heteroclinic orbit (\ref{ehetorb}) represents the figure eight 
structure as illustrated in Fig.\ref{eight}. If we denote by $S$ the 
circle, we have the topological identification:
\[
\mbox{(figure 8)}\ \otimes S = \bigcup_{p \in (-\infty,\infty),\ 
\ga \in [0,2\pi]} Q_n(p, \ga, a, \om, \pm, N)\ .
\]
\begin{figure}
\vspace{3.0in}
\caption{A geometric illustration of the figure eight structure
and its corresponding spectral identification.}
\label{eight}
\end{figure}

\newpage
\eqnsection{Coordinate-Expressions for Invariant Submanifolds}

In this section, we will give expressions for invariant submanifolds in 
coordinates, which enable us to do a Melnikov analysis and an implicit 
function argument.

\subsection{Coordinate-Expressions for Linear Invariant Submanifolds}

Consider the discrete cubic integrable nonlinear Schr{\"{o}}dinger 
equation (\ref{DNLS}), denote by $q_c$ the uniform Stokes solution,
\[
q_c = a e^{i \th(t)}\ , \ \ \ \ \th(t)=-[2(a^2-\om^2)t -\ga ] \ .
\]
Let 
\[
q_n= [a + \tq_n ] e^{i \th(t)}\ ,
\]
and linearize equation (\ref{DNLS}) at $q_c$, we have 
\begin{eqnarray*}
i \dot{\tq}_n &=& {1 \over h^2}[\tq_{n+1}-2\tq_{n}+\tq_{n-1}] \\
&+& a^2 [\tq_{n+1}+\tq_{n-1}] + 2 a^2 \overline{\tq_n}\ .
\end{eqnarray*}
Assume that $\tq_n$ takes the form,
\[
\tq_n = \bigg [ A_j e^{\Om_j t} + B_j e^{\bar{\Om}_j t} \bigg ] \cos k_j n \ ,
\]
where $k_j = 2 j \pi/ N$, ($j=0, 1, \cdot \cdot \cdot, M$). Then,
\[
\Om_j^{(\pm)} = \pm 2 \sin k_j\ \sqrt{a^2 + N^2}\sqrt{a^2-N^2
\tan^2 [k_j/2]}\ .
\]
In this paper, we only study the case that $a$ lies in the range given 
in (\ref{constr}),
so that only $\Om_1^{(\pm)}$ are real and nonzero. In fact, $\Om_0^{(\pm)}$ 
are zero, and $\Om_j^{(\pm)}$ are imaginary for ($j >1$). When $j=1$,
\[
\tq_n = a_1^{(\pm)} e^{\Om_1^{(\pm)} t} e^{\pm i \vartheta} \cos k_1 n \ ,
\]
where $a_1^{(\pm)}$ are real constants;
\begin{eqnarray*}
\vartheta &=& -{1 \over 2} \arctan \bigg \{ [ (a^2 +N^2)\sin^2(2\pi/N) \\
& & (a^2 -N^2\tan^2(\pi /N))]^{1/2}[N^2-(N^2+a^2)\cos(2\pi /N)]^{-1}
\bigg \}\ .
\end{eqnarray*}
Denote by $\B$ the block in the phase space $\cS$,
\begin{eqnarray*}
\B &\equiv& \bigg \{ \q \in \cS \ \bigg | \ q_n = e^{i\ga} \bigg [ a + 
(b_1 e^{i \vartheta} + b_2 e^{-i \vartheta}) \cos k_1 n \\
& & +\sum_{j=2}^{M} c_j \cos k_j n \bigg ]\ , \mbox{where}\ a \in 
(N\tan [\pi/N], N\tan [2\pi/N]), \ga \in [0,2\pi); \\ 
& & b_1, b_2 \ \mbox{are real}; c_j \ \mbox{is complex}; \ \mbox{and}\ 
\vartheta \ \mbox{is given above} \bigg \}\ .
\end{eqnarray*}
In terms of the coordinates $\{ a,\ga,b_1,b_2,c_j (2\leq j\leq M)\}$, 
the linear invariant center manifold $\LL^{(c)}$ is given by,
\[
\LL^{(c)} \equiv \bigg \{ \q \in \B \ \bigg | \ b_1=b_2=0 \bigg \}\ ,
\]
the linear invariant center-unstable manifold $\LL^{(cu)}$ is given by,
\[
\LL^{(cu)} \equiv \bigg \{ \q \in \B \ \bigg | \ b_2=0 \bigg \}\ ,
\]
the linear invariant center-stable manifold $\LL^{(cs)}$ is given by,
\[
\LL^{(cs)} \equiv \bigg \{ \q \in \B \ \bigg | \ b_1=0 \bigg \}\ .
\]
  
\subsection{Coordinate-Expressions for Persistent Locally 
Invariant Submanifolds}

Under the perturbed flow (\ref{PDNLS}), the linear invariant submanifolds 
$\LL^{(c)}$, $\LL^{(cu)}$ and $\LL^{(cs)}$ perturb into locally invariant 
submanifolds $W^{(c)}$, $W^{(cu)}$ and $W^{(cs)}$ (local invariance means 
that orbits can only enter or leave the submanifolds through their boundaries).
For references on the proof of such results, see for example \cite{Fen79}
\cite{LM97} \cite{LW97b}.
In terms of the coordinates $\{ a,\ga,b_1,b_2,c_j\ (2 \leq j \leq M) \}$, 
for any small $\dl_1 > 0$,
there exist $\dl_l > 0$ ($l=0,2,3$), such that $W^{(c)}$ has the expression,
\[
\left \{ \begin{array}{c} b_1 = f_1^{(c)}(a,\ga,\vc;\e,\al,\om)\ , \\ 
b_2 = f_2^{(c)}(a,\ga,\vc;\e,\al,\om)\ ; \\ \end{array} \right.
\]
$W^{(cu)}$ has the expression,
\[
b_2 = f^{(u)}(a,\ga,b_1,\vc;\e,\al,\om)\ ;
\]
and 
$W^{(cs)}$ has the expression,
\[
b_1 = f^{(s)}(a,\ga,b_2,\vc;\e,\al,\om)\ ;
\]
where $f_l^{(c)} \ (l=1,2)$, $f^{(u)}$ and $f^{(s)}$ are $C^n$ smooth 
functions for some large $n$,
$\al=(\al_1,\al_2)$, $|\e|<\dl_0$, $\ga \in[0,2\pi)$, $|b_k|<\dl_2$,
($k=1,2$), $\om \in E_{\dl_1}$, $a \in E_{\dl_1}$, $E_{\dl_1} = (
N \tan (\pi /N) +\dl_1, N \tan (2\pi /N) -\dl_1)$ when $N > 4$ and 
$=(N \tan (\pi /N) +\dl_1, \La)$ when $N=3,4$, where $\La$ is a fixed 
large constant, $\vc = (c_2, \cdot \cdot \cdot , c_M)$, 
$\| c_j \| < \dl_3\ (2 \leq j \leq M)$. 
Denote by $\A$ the annulus,
\begin{equation}
\A \equiv \bigg \{ \q \in \B \ \bigg | \ b_1=b_2=\vc=0,\ a \in E_{\dl_1}, \ 
\ga \in [0,2\pi) \bigg \}\ .
\label{anl}
\end{equation}
Then, $\A \subset \LL^{(c)}$, and $\A \subset W^{(c)}$; thus,
\[
f^{(c)}_l(a,\ga,0;\e,\al,\om)=0, \ \ (l=1,2).
\]
Notice also that $W^{(c)} = W^{(cu)} \cap W^{(cs)}$; then,
\begin{eqnarray*}
f^{(u)}(a,\ga,f^{(c)}_1(a,\ga,\vc;\e,\al,\om),\vc;\e,\al,\om)&=&f^{(c)}_2
(a,\ga,\vc;\e,\al,\om)\ , \\
f^{(s)}(a,\ga,f^{(c)}_2(a,\ga,\vc;\e,\al,\om),\vc;\e,\al,\om)&=&f^{(c)}_1
(a,\ga,\vc;\e,\al,\om)\ .
\end{eqnarray*}

\subsection{Coordinate-Expressions for Fenichel Fibers}

The center-unstable and center-stable manifolds $W^{(cu)}$ and $W^{(cs)}$ 
admit fibration through Fenichel fibers. For references on the proof of 
such results, see for example \cite{LM97} \cite{LW97b}.
Take unstable Fenichel fibers,
for example, they are a family of curves in $W^{(cu)}$, each of them is 
labeled by a base point in $W^{(c)}$. $W^{(cu)}$ is a union of the 
family of curves over $W^{(c)}$. In this sense, Fenichel fibers are 
coordinates for $W^{(cu)}$. Indeed, they are very good coordinates. 
Fenichel fibers depend $C^{n-1}$ smoothly on their base points. In backward 
time, orbits starting from points on the same fiber (in particular, from 
the base point) approach each other exponentially.

The unstable Fenichel fibers $\{ \F^{(u)}_q \}$ have the 
coordinate-expressions:
\begin{eqnarray*}
a &=& a^{(u)}(b_1;a_0,\ga_0,\vc_0;\e,\al,\om), \\
\ga &=& \ga^{(u)}(b_1;a_0,\ga_0,\vc_0;\e,\al,\om), \\
\vc &=& \vc^{\ (u)}(b_1;a_0,\ga_0,\vc_0;\e,\al,\om), \\
b_2 &=& f^{(u)}(a^{(u)},\ga^{(u)},b_1,\vc^{\ (u)};\e,\al,\om);
\end{eqnarray*}
which are $C^n$ smooth in $b_1$, and $C^{n-1}$ smooth in ($a_0,\ga_0,\vc_0$) and 
($\e,\al,\om$), where $b_1 \in (-\dl_2,\dl_2)$. The base points $q$ are 
given by,
\begin{eqnarray*}
& & a^{(u)}(f^{(c)}_1;a_0,\ga_0,\vc_0;\e,\al,\om) = a_0 \ , \\
& & \ga^{(u)}(f^{(c)}_1;a_0,\ga_0,\vc_0;\e,\al,\om)= \ga_0 \ , \\
& & \vc^{\ (u)}(f^{(c)}_1;a_0,\ga_0,\vc_0;\e,\al,\om) = \vc_0 \ , \\
& & f^{(u)}(a^{(u)},\ga^{(u)},f^{(c)}_1,\vc^{\ (u)};\e,\al,\om) \\
& & \ \ =f^{(u)}(a_0,\ga_0,f^{(c)}_1,\vc_0;\e,\al,\om) \\ 
& & \ \ =f^{(c)}_2(a_0,\ga_0,\vc_0;\e,\al,\om)\ ,
\end{eqnarray*}
where $f^{(c)}_1 = f^{(c)}_1(a_0,\ga_0,\vc_0;\e,\al,\om)$. 

The stable Fenichel fibers $\{ \F^{(s)}_q \}$ have the coordinate-expressions:
\begin{eqnarray}
a &=& a^{(s)}(b_2;a_0,\ga_0,\vc_0;\e,\al,\om), \nonumber \\
\ga &=& \ga^{(s)}(b_2;a_0,\ga_0,\vc_0;\e,\al,\om), \nonumber \\
\vc &=& \vc^{\ (s)}(b_2;a_0,\ga_0,\vc_0;\e,\al,\om), \label{sfen} \\
b_1 &=& f^{(s)}(a^{(s)},\ga^{(s)},b_2,\vc^{\ (s)};\e,\al,\om);\nonumber 
\end{eqnarray}
which are $C^n$ smooth in $b_2$, and $C^{n-1}$ smooth in ($a_0,\ga_0,\vc_0$) and 
($\e,\al,\om$), where $b_2 \in (-\dl_2,\dl_2)$. The base points $q$ are 
given by,
\begin{eqnarray*}
& & a^{(s)}(f^{(c)}_2;a_0,\ga_0,\vc_0;\e,\al,\om) = a_0 \ , \\
& & \ga^{(s)}(f^{(c)}_2;a_0,\ga_0,\vc_0;\e,\al,\om)= \ga_0 \ , \\
& & \vc^{\ (s)}(f^{(c)}_2;a_0,\ga_0,\vc_0;\e,\al,\om) = \vc_0 \ , \\
& & f^{(s)}(a^{(s)},\ga^{(s)},f^{(c)}_2,\vc^{\ (s)};\e,\al,\om) \\
& & \ \ =f^{(s)}(a_0,\ga_0,f^{(c)}_2,\vc_0;\e,\al,\om) \\ 
& & \ \ =f^{(c)}_1(a_0,\ga_0,\vc_0;\e,\al,\om)\ ,
\end{eqnarray*}
where $f^{(c)}_2 = f^{(c)}_2(a_0,\ga_0,\vc_0;\e,\al,\om)$.

\newpage
\eqnsection{Existence of Transversal Homoclinic Tubes}

In this subsection, we will establish the existence of transversal 
homoclinic tubes, 
based upon the coordinatization in the previous subsection, a Melnikov 
function calculation, and an implicit function argument.

Define the region $\U$ as follows,
\[
\U \equiv \bigg \{ \q \in \B \ \bigg | \ b_2 \in ({1 \over 2} \dl_2 , \dl_2) 
\bigg \}\ .
\] 
There exists $T > 0$ such that $F^T \circ W^{(cu)}$ intersects the region 
$\U$ (where $F^T$ is the evolution operator of the system (\ref{PDNLS})), 
the intersection $(F^T \circ W^{(cu)}) \cap \U$ has the expression,
\[
b_1 = f_T^{(u)}(a,\ga,b_2,\vc;\e,\al,\om)\ .
\]
Since for any fixed finite $T$, $F^T$ is a $C^n$ diffeomorphism, $f_T^{(u)}$ 
is also $C^n$ smooth. Notice also that when $\e =0$, $W^{(cu)}=W^{(cs)}$, 
we have
\[
f_T^{(u)}(a,\ga,b_2,\vc;0,\al,\om)=f^{(s)}(a,\ga,b_2,\vc;0,\al,\om)\ .
\]
Define the function,
\\
\begin{eqnarray*}
\tilde{\Dl}(a_0,\ga_0,\vc_0,b_2;\e,\al,\om) &=& f^{(u)}_T(a^{(s)},\ga^{(s)},
b_2,\vc^{\ (s)};\e,\al,\om) \\ 
&-&  f^{(s)}(a^{(s)},\ga^{(s)},b_2,\vc^{\ (s)};\e,\al,\om) \ ,
\end{eqnarray*}
\\
where
\begin{eqnarray*}
a^{(s)} &=& a^{(s)}(b_2;a_0,\ga_0,\vc_0;\e,\al,\om), \\
\ga^{(s)} &=& \ga^{(s)}(b_2;a_0,\ga_0,\vc_0;\e,\al,\om), \\
\vc^{\ (s)} &=& \vc^{\ (s)}(b_2;a_0,\ga_0,\vc_0;\e,\al,\om);
\end{eqnarray*}
are defined in the coordinate-expression (\ref{sfen}) for the stable 
Fenichel fibers, and $\tilde{\Dl}$ is a $C^{n-1}$ smooth function. 
Setting $b_2 = {3 \over 4} \dl_2$, for sufficiently small 
$\e$ and $\| \vc_0 \ \|$ ($\| \vc_0 \ \|=[\sum_{j=2}^M|c_j|^2]^{1/2}$), 
we have 
\\
\begin{eqnarray*}
\tilde{\Dl}(a_0,\ga_0,\vc_0,{3 \over 4} \dl_2;\e,\al,\om) &=& 
\e M_{\tF_1}(a_0,\ga_0;\al,\om) +\e R_1(a_0,\ga_0,\vc_0;\e,\al,\om)  \\
&+& \e^2 R_2(a_0,\ga_0,\vc_0;\e,\al,\om) \ ,
\end{eqnarray*}
where $M_{\tF_1}(a_0,\ga_0;\al,\om)$ is the Melnikov function,
\[
R_1 \sim O(\| \vc_0\|)\ , \ \ \mbox{as}\ \| \vc_0 \| \ra 0\ .
\]
For a derivation on this result, see for example \cite{LM97}.
Define the function $\Dl$ as follows,
\\
\begin{eqnarray}
\Dl &=& \Dl(a_0,\ga_0,\vc_0;\e,\al,\om) \nonumber \\
    &=& {1 \over \e} \tilde{\Dl}(a_0,\ga_0,\vc_0,{3 \over 4} \dl_2;\e,\al,\om)
\nonumber \\
&=& M_{\tF_1}(a_0,\ga_0;\al,\om) + R_1(a_0,\ga_0,\vc_0;\e,\al,\om) 
\label{melm} \\
&+& \e R_2(a_0,\ga_0,\vc_0;\e,\al,\om)\ .\nonumber
\end{eqnarray}
Then $\Dl$ is $C^{n-2}$ smooth in ($a_0,\ga_0,\vc_0;\e,\al,\om$). The 
Melnikov function $M_{\tF_1}=M_{\tF_1}(a,\ga;\al,\om)$ is given as:
\[
M_{\tF_1}=i \int_{-\infty}^{\infty}\{ \tF_1, H_1 \}|_{\vQ} dt=2
\int_{-\infty}^{\infty}\sum_{n=0}^{N-1}\ \mbox{Im}\ \bigg \{ 
{\pa \tF_1 \over \pa r_n} \rho_n {\pa H_1 \over \pa q_n}\bigg \}
_{\vQ} dt\ ,
\]
where 
\[
\{ \tF_1, H_1 \}=\sum_{n=0}^{N-1}\bigg [{\pa \tF_1 \over \pa q_n} \rho_n 
{\pa H_1 \over \pa r_n} -{\pa \tF_1 \over \pa r_n} \rho_n 
{\pa H_1 \over \pa q_n}\bigg ]\ ,
\]
and $\vQ$ is given in (\ref{ehetorb}) and $\mbox{grad}\ \tF_1(\vQ)$
is given in (\ref{emelv}). Thus,
\begin{eqnarray}
{1 \over 2} M_{\tF_1} &=& \al_1 \int_{-\infty}^{\infty}
\sum_{n=0}^{N-1}\ \mbox{Im}\ \bigg \{ {\pa \tF_1 \over \pa r_n}
\bigg [ {\rho_n \over h^2}\ln \rho_n +(q_n +\bq_n)\bq_n \bigg ]
\bigg \}_{\vQ} dt \nonumber \\
&+& \al_2 \int_{-\infty}^{\infty}
\sum_{n=0}^{N-1}\ \mbox{Im}\ \bigg \{ {\pa \tF_1 \over \pa r_n}
\bigg [ 2 q_n {\rho_n \over h^2}\ln \rho_n +(q_n^2+\bq_n^2)\bq_n\bigg ]
\bigg \}_{\vQ} dt \label{melf} \\
&\equiv& f_1(a,\ga;\om,N) \al_1 + f_2(a,\ga;\om,N) \al_2\ .\nonumber
\end{eqnarray}
Setting  $M_{\tF_1}=0$ in (\ref{melf}), we have an equation of the form,
\begin{equation}
\al_1 - 4 \om \k \al_2 = 0 \ , \label{mrt}
\end{equation}
where
\[
\k = \k(a,\ga;\om,N) = -{1 \over 4 \om}{f_2(a,\ga;\om,N) \over 
f_1(a,\ga;\om,N)}\ .
\]
The graphs of $\k$ are shown in figures \ref{kap1} and \ref{kap2}. 
Figure \ref{kap2} is for the case when $q_c$ is a circle of fixed points.
\begin{figure}
\vspace{1.5in}
\caption{The graph of $\k$ for a fixed value of $\om$.}
\label{kap1}
\end{figure}
\begin{figure}
\vspace{1.5in}
\caption{The graph of $\k$ when $a=\om$.}
\label{kap2}
\end{figure}
For the convenience of later argument, we denote $M_{\tF_1}$ simply 
by $M$ and denote the surface defined by (\ref{mrt}) by $S_\ga$:
\begin{equation}
S_\ga \ \ \ :\ \ \ \ga_0 = \Ga^{(0)}(a_0;\al,\om)\ .
\label{srf}
\end{equation}
\begin{theorem}
There exist a positive constant $\e_0 > 0$ and a region $\E$ for ($\al,\om$), 
such that for any $\e \in (-\e_0,\e_0)$ and any $(\al,\om) \in \E$, there 
exists a codimension 2 transversal homoclinic tube asymptotic to 
the codimension 2 center manifold $W^{(c)}$.
\label{thmht}
\end{theorem}

Proof: There exists a region $\V$ on the surface $S_\ga$ in (\ref{srf}), such 
that
\[
{\pa \over \pa \ga_0}M(a_0,\ga_0;\al,\om) \neq 0
\]
and is bounded. See figures (\ref{der1}) and (\ref{der2}) for the 
corresponding graphs of ${\pa M \over \pa \ga_0}$. 
\begin{figure}
\vspace{1.5in}
\caption{The graph of ${\pa M \over \pa \ga_0}$ for a fixed value of $\om$.}
\label{der1}
\end{figure}
\begin{figure}
\vspace{1.5in}
\caption{The graph of ${\pa M \over \pa \ga_0}$ when $a = \om$.}
\label{der2}
\end{figure}
Next we want to solve 
the equation (\ref{melm}) by the implicit function theorem [\cite{Die60},p265].
For any $(a^*_0,\ga^*_0;\al^*,\om^*) \in \V$,
\begin{equation}
{\pa \over \pa \ga_0}\Dl(a^*_0,\ga^*_0,0;0,\al^*,\om^*) 
= {\pa \over \pa \ga_0}M(a^*_0,\ga^*_0;\al^*,\om^*) \neq 0
\label{par}
\end{equation}
and is bounded. Then by the implicit function theorem 
[\cite{Die60},p265], there is a neighborhood $\W^{(*)}$ of ($a_0^*, \vc_0=0; 
\e =0, \al^*, \om^*$) and a unique $C^{n-2}$ function,
\[
\ga_0 = \Ga^{(*)}(a_0,\vc_0;\e,\al,\om)
\]
defined in $\W^{(*)}$, such that 
\[
\Ga^{(*)}(a^*_0,\vc_0=0;\e=0,\al^*,\om^*) = \ga_0^* \ ,
\]
and 
\[
\Dl(a_0,\Ga^{(*)}(a_0,\vc_0;\e,\al,\om),\vc_0;\e,\al,\om) = 0 \ .
\]
Since $\Dl$ is a $C^{n-2}$ smooth function, by relation (\ref{par}), we have 
\[
{\pa \over \pa \ga_0}\Dl(a_0,\Ga^{(*)}(a_0,\vc_0;\e,\al,\om),
\vc_0;\e,\al,\om)\neq 0\ ,
\]
and is bounded for $(a_0,\vc_0;\e,\al,\om) \in \W^{(*)}$. Thus, the 
center-unstable 
manifold $W^{(cu)}$ and the center-stable manifold $W^{(cs)}$ have a 
transversal intersection at the neighborhood $\W^{(*)}$.
Let 
\[
\W = \bigcup_{(a^*_0,\ga^*_0;\al^*,\om^*) \in \V} \W^{(*)} \ ;
\]
then there is a unique $C^{n-2}$ function,
\begin{equation}
\ga_0 = \Ga(a_0,\vc_0;\e,\al,\om)
\label{sface}
\end{equation}
defined in $\W$, such that 
\[
\Ga(a_0,\vc_0=0;\e=0,\al,\om) = \Ga^{(0)}(a_0;\al,\om) \ ,
\]
and 
\[
\Dl(a_0,\Ga(a_0,\vc_0;\e,\al,\om),\vc_0;\e,\al,\om) = 0 \ .
\]
Notice that (\ref{sface}) defines a codimension 1 submanifold $W^{(c)}_b$ of 
the center manifold $W^{(c)}$, which has the expression,
\begin{eqnarray*}
b_1 &=& f_1^{(c)}(a_0,\Ga(a_0,\vc_0;\e,\al,\om),\vc_0;\e,\al,\om)\ , \\
b_2 &=& f_2^{(c)}(a_0,\Ga(a_0,\vc_0;\e,\al,\om),\vc_0;\e,\al,\om)\ , \\
\ga_0 &=& \Ga(a_0,\vc_0;\e,\al,\om)\ ; 
\end{eqnarray*}
where $(a_0,\vc_0;\e,\al,\om) \in \W$. Define $\HH$ as follows,
\[
\HH = \bigcup_{t \in (-\infty,\infty)} F^t \circ 
\bigcup_{q \in W^{(c)}_b} \F^{(s)}_q \ ,
\]
where $F^t$ is the evolution operator for the system (\ref{PDNLS}).
Then $\HH$ is the codimension 2 transversal homoclinic tube asymptotic 
to $W^{(c)}$. $\Box$

For an illustration of the homoclinic tube, see figure \ref{fight}. Studies 
on the symbolic dynamics in the neighborhood of this homoclinic tube 
are topics of future works.
\begin{figure}
\vspace{1.5in}
\caption{An illustration of the homoclinic tube.}
\label{fight}
\end{figure}

\newpage
\eqnsection{Symbolic Dynamics of Segments for the Case ($N=3$, Nonresonant 
Region)}

Denote by $\Pi$ the plane,
\begin{equation}
\Pi \equiv \bigg \{ \q \in \cS \ \bigg | \ q_n = q\ \forall n \bigg \} \ .
\label{tubepi}
\end{equation}
$\Pi$ is invariant under the flow governed by (\ref{PDNLS}). 
The annulus $\A$ defined in (\ref{anl}) is a subset of $\Pi$. Away from 
the resonant circle:
\begin{equation}
S_\ga \equiv \bigg \{ \q \in \B \ \bigg | \ b_1=b_2=\vc=0,\ a = \om, \ 
\ga \in [0,2\pi) \bigg \}\ ,
\label{cir}
\end{equation}
the dynamics is given by periodic motions, and denote by $\hat{\A}$ 
such an annulus with boundaries which are periodic orbits. When $N=3$, 
the Melnikov analysis shows that there exists a transversal homoclinic 
tube asymptotic to $\Pi$. By the invariance of the Hamiltonian $H$,
each orbit inside the homoclinic tube, approaching a periodic orbit in 
either the forward or the backward time, has to approach the same periodic 
orbit in both forward and backward time. Thus, there exists a homoclinic 
tube asymptotic to the annulus $\hat{\A}$, as illustrated in figure 
\ref{n3ht}.
\begin{definition}
Define a Poincar\'{e} section $\Sg$ as follows,
\[
\Sg = \bigg \{ \q \in \B \ \bigg | \ \ga = 0 \bigg \} \ ,
\]
and let $P$ be the Poincar\'{e} map,
\[
P\ :\ \D \subset \Sg \mapsto \Sg \ ,
\]
induced by the flow, where $\D$ is the domain of definition for $P$.
\end{definition}
The transversal homoclinic tube $\HH_P$ for the Poincar\'{e} map $P$ 
is illustrated in figure \ref{n3pht}.
\begin{figure}
\vspace{1.5in}
\caption{A homoclinic tube asymptotic to a nonresonant annulus when 
$N=3$.}
\label{n3ht}
\end{figure}
\begin{figure}
\vspace{1.5in}
\caption{A homoclinic tube asymptotic to a segment under the 
Poincar\'{e} map $P$ when $N=3$.}
\label{n3pht}
\end{figure}
This transversal homoclinic tube $\HH_P$ consists of segments, 
and is asymptotic to a segment $s$ in the annulus $\hat{\A}$. 
The segment $s$ consists of fixed points of $P$. Invariant tubes 
in the neighborhood of $\HH_P$ also consist of segments. We have 
the following symbolic dynamics theorem for the homoclinic tube 
$\HH_P$.
\begin{theorem}[Silnikov \cite{Sil68}]
When $N=3$,
the set of invariant tubes of the Poincar\'{e} map $P$ lying 
wholly within a sufficiently small neighborhood of the transversal 
homoclinic tube $\HH_P$ is in one-to-one correspondence with 
the set of all doubly infinite sequences,
\[
J=(\cdot \cdot \cdot, j_{-1}, j_0, j_1, \cdot \cdot \cdot)\ ,
\]
where $j_l \geq j_*$ for all $l$ and some large integer $j_*$.
\end{theorem}

\newpage
\eqnsection{Structures Inside the Asymptotic Manifolds of the Transversal 
Homoclinic Tubes}

Structures inside the transversal homoclinic tubes can be very complicated. 
Since we are studying near integrable Hamiltonian systems, the asymptotic 
manifold (i.e. the center manifold) $W^{(c)}$ often contains both KAM tori 
and stochastic layers. Studies on such structures are topics of future 
works. In this section, we study interesting structures on an invariant 
plane inside the center manifold $W^{(c)}$, and special orbits inside 
the transversal homoclinic tubes.

The dynamics on $\Pi$ (\ref{tubepi}) is governed by 
\begin{equation}
i \dot{q} = 2 [|q|^2 -\om^2] q + \e \bigg \{ [\al_1 (q +\bq) +\al_2 (q^2 
+\bq^2)]q +[\al_1 +2\al_2 \bq]{\rho \over h^2} \ln \rho \bigg \}\ , 
\label{CPNLS}
\end{equation}
where $\rho = 1+h^2|q|^2$.
Let $q = I e^{i \xi}$, where $I$ is the modulus of $q$; then,
\begin{eqnarray*}
dI/dt &=& - \e \sin \xi \ [\al_1 +4 \al_2 I \cos \xi ] {\rho \over h^2}
\ln \rho \ , \\
d\xi/dt &=& -2[I^2-\om^2]-\e \bigg \{ 2 \al_1 I\cos \xi +2\al_2 I^2 \cos 2\xi
\\
& & +\bigg [\al_1{\cos \xi \over I} +2\al_2 \cos 2 \xi \bigg ]{\rho \over h^2}
\ln \rho \bigg \}\ .
\end{eqnarray*}
The dynamics away from the neighborhood of $I=\om$ is given by 
periodic motions.
The dynamics in the neighborhood of $I=\om$ is very interesting. Let 
$I=\om +\eta y$, $\eta = \sqrt{\e}$, and define the 
{\em{resonant annulus}} $\tilde{\A}$,
\[
\tilde{\A} \equiv \bigg \{ (y,\xi)\ \bigg | \ \xi \in [0,2\pi], y \ 
\mbox{is bounded} \bigg \}.
\]
Inside the resonant annulus $\tilde{\A}$, the dynamics is governed by,
\begin{eqnarray}
dy/d\tau &=& f_1\ , \nonumber \\
\label{rcdy} \\
d\xi/d\tau &=& f_2 \ , \nonumber
\end{eqnarray}
where
\begin{eqnarray}
f_1 &=& -\sin \xi \ [\al_1 +4 \al_2 (\om + \eta y) \cos \xi ]
{\rho \over h^2}\ln \rho\ , \nonumber \\
f_2 &=& -4 \om y - \eta \bigg \{ 2y^2  + 2 \al_1 (\om +\eta y) 
\cos \xi \nonumber \\
& & + 2 \al_2 \cos 2 \xi \ (\om +\eta y)^2 +\bigg [ {\al_1 \cos \xi 
\over \om +\eta y} +2\al_2 \cos 2\xi \bigg ] {\rho \over h^2}\ln \rho 
\bigg \}\ ,\nonumber
\end{eqnarray}
where $\tau = \eta t$.

\subsection{First Order Dynamics in the Resonant Annulus}

Setting $\eta=0$ in (\ref{rcdy}), we have the system,
\begin{eqnarray}
dy/d\tau &=& -\Om \sin \xi \ [\al_1 +4 \al_2 \om \cos \xi ]\ , 
\nonumber \\
\label{urcdy} \\
d\xi/d\tau &=& -4 \om y \ ,\nonumber 
\end{eqnarray}
where $\Om = h^{-2}\rho_0 \ln \rho_0$ and $\rho_0 =1+h^2\om^2$.
This first order system (\ref{urcdy}) is also a Hamiltonian system,
\begin{eqnarray*}
dy/d\tau &=& {\pa \tilde{H} \over \pa \xi}\ , \\
d\xi/d\tau &=& -{\pa \tilde{H} \over \pa y}\ ; 
\end{eqnarray*}
where $\tilde{H}=2 \om y^2 +\Om (\al_1 \cos \xi +\al_2 \om \cos 2 \xi)$.
The fixed points of the system (\ref{urcdy}) are as follows:
\begin{itemize}
\item When $\bigg | {\al_1 \over 4 \al_2 \om } \bigg | < 1$, there are 
four fixed points,
\begin{eqnarray}
& & y^{(j)}_0 = 0 \ , \ \ (j=1,2,3,4)\ ; \nonumber \\
\label{ufpt1} \\
& & \xi^{(1)}_0 = 0 \ , \ \ \xi^{(2)}_0 = \pi\ , \ \ \xi^{(3)}_0 = -
\xi^{(4)}_0=\arccos \bigg [ -{\al_1 \over 4 \al_2 \om } \bigg ] \ .
\nonumber 
\end{eqnarray}
\item When $\bigg | {\al_1 \over 4 \al_2 \om } \bigg | > 1$, there are 
two fixed points,
\begin{eqnarray}
& & y^{(j)}_0 = 0 \ , \ \ (j=1,2)\ ; \nonumber \\
\label{ufpt2} \\
& & \xi^{(1)}_0 = 0 \ , \ \ \xi^{(2)}_0 = \pi\ .\nonumber 
\end{eqnarray}
\end{itemize}
The phase diagrams are shown in figure \ref{4pd} when $\bigg | 
{\al_1 \over 4 \al_2 \om } \bigg | < 1$, and in figure \ref{2pd} 
when $\bigg | {\al_1 \over 4 \al_2 \om } \bigg | > 1$.
\begin{figure}
\vspace{3.0in}
\caption{Phase diagrams in the resonant annulus when $\bigg | 
{\al_1 \over 4 \al_2 \om } \bigg | < 1$, (a). $\al_2 >0$, $\al_1 >0$;
(b). $\al_2 >0$, $\al_1 <0$; (c). $\al_2 <0$, $\al_1 >0$;
(d). $\al_2 <0$, $\al_1 <0$.}
\label{4pd}
\end{figure}
\begin{figure}
\vspace{1.5in}
\caption{Phase diagrams in the resonant annulus when $\bigg | 
{\al_1 \over 4 \al_2 \om } \bigg | > 1$, (a). $\al_1 >0$;
(b). $\al_1 <0$.}
\label{2pd}
\end{figure}

\subsection{The Full Dynamics in the Resonant Annulus}

First we prove the existence of fixed points for the system (\ref{rcdy})
governing the full dynamics in the resonant annulus. 
\begin{lemma}
For any large $\La > 0$, any small $\dl_0 > 0$, there exists 
$\eta_0 > 0$, such that there exist fixed points for the 
system (\ref{rcdy}),
\[
y^{(j)}= y^{(j)}(\eta,\al_1,\al_2,\om)\ , \ \ \ \ 
\xi^{(j)}= \xi^{(j)}(\eta,\al_1,\al_2,\om)\ ;
\]
which are $C^1$ smooth, where ($j=1,2,3,4;$ when $\bigg |{\al_1 
\over 4\om \al_2}\bigg | < 1-\dl_0$) and ($j=1,2;$ when 
$\bigg |{\al_1 \over 4\om \al_2}\bigg | > 1+\dl_0$), 
$|\eta| < \eta_0$, $\dl_0 < |\al_j| < \La $, ($j=1,2$), 
$\dl_0 < |\om|< \La$; and 
\begin{eqnarray*}
& & y^{(j)}(0,\al_1,\al_2,\om) =y^{(j)}_0=y^{(j)}_0(\al_1,\al_2,\om)\ , \\  
& & \xi^{(j)}(0,\al_1,\al_2,\om)=\xi^{(j)}_0=\xi^{(j)}_0(\al_1,\al_2,\om)\ .
\end{eqnarray*}
\label{fiptl}
\end{lemma}

Proof: We want to solve the equations
\[
f_1(y,\xi;\eta,\al_1,\al_2,\om)= 0\ , \ \ 
f_2(y,\xi;\eta,\al_1,\al_2,\om)= 0\ ;
\]
in the neighborhoods of the fixed points
\begin{equation}
X^{(j)}_0=(y^{(j)}_0,\xi^{(j)}_0;0,\al_1,\al_2,\om)\ .
\label{ufipt}
\end{equation}
Notice that
\begin{eqnarray*}
& & {\pa f_1 \over \pa y}(y,\xi;0,\al_1,\al_2,\om) = 0 \ ,\\
& & {\pa f_1 \over \pa \xi}(y,\xi;0,\al_1,\al_2,\om) \\
& & ={\rho_0 \over h^2} \ln \rho_0 \ \bigg [ 4\al_2 \om \sin^2 \xi 
-\cos \xi \ (\al_1 +4 \om \al_2 \cos \xi) \bigg ]\ , \\
& & {\pa f_2 \over \pa y}(y,\xi;0,\al_1,\al_2,\om) = -4\om \ , \\
& & {\pa f_2 \over \pa \xi}(y,\xi;0,\al_1,\al_2,\om) = 0\ .
\end{eqnarray*}
Then for any $\dl_0 > 0$,
\begin{eqnarray*}
& & \ \mbox{det}\ \left | \begin{array}{lr} {\pa f_1 \over \pa y}(X^{(j)}_0) & 
{\pa f_1 \over \pa \xi}(X^{(j)}_0) \cr \cr {\pa f_2 \over \pa y}(X^{(j)}_0) & 
{\pa f_2 \over \pa \xi}(X^{(j)}_0) \cr \end{array} \right | \\
&=& {4\om \rho_0 \over h^2}\ln \rho_0\  \bigg [ 4 \al_2 \om 
\sin^2 \xi_0^{(j)} - 
\cos \xi_0^{(j)}\ (\al_1 +4\al_2 \om \cos \xi_0^{(j)}) \bigg ] \neq 0 \ ,
\end{eqnarray*}
when $\bigg |{\al_1 \over 4\om \al_2}\bigg | < 1-\dl_0$ or 
$\bigg |{\al_1 \over 4\om \al_2}\bigg | > 1+\dl_0$. Then
\[
\left ( \begin{array}{lr} {\pa f_1 \over \pa y}(X^{(j)}_0) & 
{\pa f_1 \over \pa \xi}(X^{(j)}_0) \cr \cr {\pa f_2 \over \pa y}(X^{(j)}_0) & 
{\pa f_2 \over \pa \xi}(X^{(j)}_0) \cr \end{array} \right )
\]
define linear homeomorphisms. Notice that $f_1$ and $f_2$ are 
$C^1$ functions;
then by the implicit function theorem [\cite{Die60}, pp265], for 
any fixed value of $X^{(j)}_0$ which satisfies the restrictions 
given in the Lemma, there is a neighborhood $V^{(j)}_Y$ of 
$Y=(\eta = 0, \al_1,\al_2,\om)$, and unique functions
\begin{equation}
y^{(j)} = y^{(j)}(\eta, \al_1,\al_2,\om), \ \ \xi^{(j)} = 
\xi^{(j)}(\eta, \al_1,\al_2,\om);
\label{rpq3}
\end{equation}
which are $C^1$ functions defined in $V^{(j)}_Y$, such that 
\[
f_l(y^{(j)},\xi^{(j)};\eta,\al_1,\al_2,\om)=0, \ (l=1,2)\ .
\]
Let $V^{(j)}=\bigcup_Y V^{(j)}_Y$; then the functions (\ref{rpq3}) 
are uniquely defined in $V^{(j)}$. There exists $\eta_0 > 0$, such 
that $|\eta| < \eta_0$ and the rest of restrictions in the Lemma 
define subregions of $V^{(j)}$. $\Box$

In fact, we have
\[
\xi^{(1)}=0\ , \ \ \ \ \xi^{(2)}=\pi\ ,
\]
and $y^{(j)}\ (j=1,2)$ satisfy the equations,
\[
I^2 - \om^2 +\eta^2 \bigg [ \al_2 I^2 \pm \al_1 I +\bigg (\al_2 \pm 
{\al_1 \over 2I}\bigg ) {\rho \over h^2} \ln \rho \bigg ] = 0\ ,
\]
where $I=\om +\eta y$, `$+$' for $j=1$ and `$-$' for $j=2$;
\[
\xi^{(3)}=-\xi^{(4)}=\arccos \bigg [ - {\al_1 \over 4\al_2 I} \bigg ]\ ,
\]
and $y^{(3)}=y^{(4)}$ satisfies the equation,
\[
(1-\eta^2 \al_2)I^2 - (\om^2 +{\eta^2 \al_1^2 \over 8 \al_2}) - \eta^2 
\al_2  {\rho \over h^2} \ln \rho = 0\ ,
\]
where $I=\om +\eta y$. Approximate expressions of these fixed points can be 
obtained:
\begin{equation}
y^{(j)} = \eta y^{(j)}_1 + O(\eta^2)\ , \ \ (j=1,2); 
\label{rpel3}
\end{equation}
where 
\begin{eqnarray*}
y^{(1)}_1 &=& - {1 \over 4 \om} \bigg [ 2 \om (\al_2 \om +\al_1) +(2 \al_2 
+\al_1/\om) {\rho_0 \over h^2} \ln \rho_0 \bigg ] \ , \\
y^{(2)}_1 &=& - {1 \over 4 \om} \bigg [ 2 \om (\al_2 \om -\al_1) +(2 \al_2 
-\al_1/\om) {\rho_0 \over h^2} \ln \rho_0\bigg ] \ .
\end{eqnarray*}
When $\bigg |{\al_1 \over 4\om \al_2}\bigg | < 1$,
\begin{equation}
y^{(l)} = \eta y^{(l)}_1 + O(\eta^2)\ , \ \ (l=3,4); 
\label{rpel4}
\end{equation}
where
\[
y^{(3)}_1 = y^{(4)}_1 ={1 \over 16 \al_2 \om}(8\al_2^2 \om^2 +\al_1^2)
+{\al_2 \over 2 \om} {\rho_0 \over h^2}\ln \rho_0\ . 
\]
\begin{eqnarray*}
\xi^{(3)}=-\xi^{(4)}&=& \arccos \bigg [ -{\al_1 \over 4\al_2 \om} 
+\eta^2 \psi +\ O(\eta^3) \bigg ] \\
&=& \xi^{(3)}_0 - \eta^2 \bigg [ 1-{\al_1^2 \over 16 \al_2^2 \om^2}
\bigg ]^{-1/2} \psi +\ O(\eta^3)\ ,
\end{eqnarray*}
where $\psi = {\al_1 \over 4\al_2 \om^2}y_1^{(3)}$.
Let $X^{(j)} = (y^{(j)},\xi^{(j)};\eta,\al_1,\al_2,\om)$ denote 
these fixed points and notice that $X^{(j)}_0 = (y^{(j)}_0,
\xi^{(j)}_0;\eta=0,\al_1,\al_2,\om)$ a notation used before in (\ref{ufipt}).
\begin{lemma}
The fixed points $X^{(j)}$ have the same saddle or center nature 
as $X^{(j)}_0$
\label{prot}
\end{lemma}

Proof: 
\begin{eqnarray*}
{\pa f_1 \over \pa y} &=& -\sin \xi \ [4\al_2 \sqrt{\e} \cos \xi ]{\rho \over 
h^2}\ln \rho \\
& & -\sin \xi \ [\al_1 +4\al_2 (\om +\sqrt{\e} y)\cos \xi ]{d \over dy}\bigg ( 
{\rho \over h^2}\ln \rho \bigg )\ , \\
{\pa f_2 \over \pa \xi} &=& \sin \xi \ [4\al_2 \sqrt{\e} \cos \xi ]{\rho \over 
h^2}\ln \rho \\
& & +\sqrt{\e}\sin \xi \ [\al_1 +4\al_2 (\om +\sqrt{\e} y)\cos \xi ]  \\
& & \bigg \{ 2 (\om +\sqrt{\e} y) +{1 \over \om +\sqrt{\e} y}{\rho \over 
h^2}\ln \rho \bigg \}\ .
\end{eqnarray*}
Then,
\begin{equation}
{\pa f_1 \over \pa y}(X^{(j)}) +{\pa f_2 \over \pa \xi}(X^{(j)}) = 0 \ ,
\label{rpel5}
\end{equation}
in fact, when $j=1,2$, ${\pa f_1 \over \pa y}(X^{(j)}) 
={\pa f_2 \over \pa \xi}(X^{(j)}) = 0$. 
Linearize (\ref{rcdy}) at $X^{(j)}$, we have 
\begin{equation}
{d \over d\tau}\left ( \begin{array}{c} \tilde{y} \\ \\ \tilde{\xi} \\ 
\end{array} \right )
=L_j \ \left ( \begin{array}{c} \tilde{y} \\ \\ \tilde{\xi} \\ \end{array} 
\right )\ ,
\label{rpel6}
\end{equation}
where
\[
L_j =\left ( \begin{array}{lr} \frac {\pa f_1}{\pa y}(X^{(j)}) &  
\frac {\pa f_1}{\pa \xi}(X^{(j)}) \\ \\ \frac {\pa f_2}{\pa y}(X^{(j)}) &  
\frac {\pa f_2}{\pa \xi}(X^{(j)}) \\ \end{array} \right )\ ,
\]
where
\begin{eqnarray}
{\pa f_l \over \pa y}(X^{(j)}) &=& {\pa f_l \over \pa y}(X^{(j)}_0)+ O(\eta)\ ,
\nonumber \\
& &  \ \ \ \ \ \ \ \ \ \ \ \ \ \ \ \ \ \ \ \ \ \ \ \ \ (l=1,2) \ \ 
\label{rpel7} \\ 
{\pa f_l \over \pa \xi}(X^{(j)}) &=& {\pa f_l \over \pa \xi}(X^{(j)}_0)+ 
O(\eta)\ . \nonumber
\end{eqnarray}
The eigenvalues of $L_j$ satisfy 
\[
\la^2 - \ \mbox{trace}(L_j) \la + \ \mbox{det}(L_j) = 0.
\]
By relations (\ref{rpel5}) and (\ref{rpel7}),
\begin{eqnarray*}
\mbox{trace}(L_j) &=& 0 \ , \\ 
\mbox{det}(L_j) &=& \ \mbox{det}\left ( \begin{array}{lr} 
\frac {\pa f_1}{\pa y}(X^{(j)}_0) &  
\frac {\pa f_1}{\pa \xi}(X^{(j)}_0) \\ \\ \frac{\pa f_2}{\pa y}(X^{(j)}_0) &  
\frac {\pa f_2}{\pa \xi}(X^{(j)}_0) \\ \end{array} \right ) + O(\eta)\ ;
\end{eqnarray*}
thus, when $\eta$ is sufficiently small, the fixed points $X^{(j)}$ have 
the same types of stability as $X^{(j)}_0$. $\Box$

The phase diagram is given by the level sets of the rescaled Hamiltonian,
\[
\hat{H} = {\rho_0 \over 2N\om \eta^2}\bigg [ H - {2N \over h^2} (\om^2 
-{\rho_0 \over h^2} \ln \rho_0) \bigg ]\ ,
\]
where $H$ is the restriction of the Hamiltonian to the resonant annulus 
$\tilde{\A}$,
\begin{eqnarray*}
H &=& {2N \over h^2} \bigg \{ [(\om +\eta y)^2-{\rho_0 \over h^2}\ln \rho ] \\
& &+\eta^2[\al_1 (\om +\eta y)\cos \xi +\al_2(\om +\eta y)^2\cos 2\xi ]
\ln \rho \bigg \}\ .
\end{eqnarray*}
$\hat{H}$ is smooth in $\eta$ and has the approximate expression,
\[
\hat{H}=2\om y^2 +[\al_1 \cos \xi +\al_2 \om \cos 2\xi ]{\rho_0 \over h^2}
\ln \rho_0 +\ O(\eta)\ .
\]
Setting $\eta=0$, we have the phase diagrams as shown in 
figures \ref{4pd} and \ref{2pd}. When $\eta \neq 0$, invariant manifolds 
of the saddles perturb smoothly, and are given by level sets of $\hat{H}$; 
thus, the figure eight loops do not break. By Lemmas \ref{fiptl} and 
\ref{prot}, in the 
resonant annulus $\tilde{\A}$, the phase diagrams when $\eta \neq 0$ are 
topologically equivalent to those when $\eta = 0$.

\subsection{Special Orbits Inside the Transversal Homoclinic Tubes}

As a consequence of Theorem \ref{thmht}, there exist orbits asymptotic 
to fixed points in the resonant annulus $\tilde{\A}$ in forward or 
backward time. See figures \ref{sporb1} and \ref{sporb2} for some examples.
Whether or not these orbits are actually homoclinic orbits or heteroclinic 
orbits is a very difficult open question.
\begin{figure}
\vspace{3.0in}
\caption{Special orbits asymptotic to fixed points in the resonant annulus 
in forward or backward time when $\bigg | 
{\al_1 \over 4 \al_2 \om } \bigg | < 1$, (a). $\al_2 >0$, $\al_1 >0$;
(b). $\al_2 >0$, $\al_1 <0$; (c). $\al_2 <0$, $\al_1 >0$;
(d). $\al_2 <0$, $\al_1 <0$.}
\label{sporb1}
\end{figure}
\begin{figure}
\vspace{1.5in}
\caption{Special orbits asymptotic to fixed points in the resonant annulus 
in forward or backward time when $\bigg | 
{\al_1 \over 4 \al_2 \om } \bigg | > 1$, (a). $\al_1 >0$;
(b). $\al_1 <0$.}
\label{sporb2}
\end{figure}

\newpage
\subsection{Homoclinic and Heteroclinic Orbits for the 
Case ($N=3$, Resonant Region)}

Inside the resonant annulus $\tilde{\A}$, the dynamics for both the full 
system (\ref{rcdy}) and the first order system (\ref{urcdy}) is shown in 
figures \ref{4pd} and \ref{2pd}. As a consequence of Theorem \ref{thmht}, 
there exist orbits asymptotic to the saddles in forward time. For fixed 
($\e,\al,\om$), these orbits correspond to the intersection points 
between the surface (\ref{sface}) and the figure eight level sets 
of the saddles. In backward time, these orbits approach the resonant 
annulus $\tilde{\A}$. Since the Hamiltonian is conserved and the figure 
eights are given as the level sets of the Hamiltonian, in backward time, 
these orbits also approach the saddles on the same figure eight level sets. 
Thereby, we have homoclinic or heteroclinic orbits as shown in 
figure \ref{hhorb} for example.
\begin{figure}
\vspace{1.5in}
\caption{Homoclinic and heteroclinic orbits for the 
case ($N=3$, resonant region).}
\label{hhorb}
\end{figure}
In drawing these figures, the fact that the quantity
\[
I_2=\sum_{n=0}^{N-1} [\bq_n (q_{n+1}+q_{n-1})]
\]
only changes $O(\e)$ in finite time interval, and the fibers are 
smooth with respect to their base points and $\e$.

\newpage
\eqnsection{Conclusion}

In this paper, we have proved the existence of homoclinic tubes 
for the discrete cubic nonlinear Schr\"odinger equation under Hamiltonian 
perturbations, and discuss the symbolic dynamics of invariant tubes 
in the neighborhoods of such homoclinic tubes when the system is four 
dimensional. In future works, we will study such symbolic dyanmics 
structures when the dimension is higher. Since the system 
that we study here is a near-integrable Hamiltonian system,
we are also interested in studying KAM 
tori and stochastic layer structures inside each invariant tube, 
in particular, the homoclinic tube, and the connection with the 
symbolic dynamics of invariant tubes. Due to the normally hyperbolic
nature of the center manifold, the relevant work on such KAM tori 
is that of Graff \cite{Gra74}.

{\bf Acknowlegment}: I am greatly indebted to Brenda Frazier for 
compiling all the artistic works.

\newpage
\bibliography{TDNLS}

\begin{thebibliography}{10}

\bibitem{AL76}
M.~J. Ablowitz and J.~F. Ladik.
\newblock A {N}onlinear {D}ifference {S}cheme and {I}nverse {S}cattering.
\newblock {\em Stud. Appl. Math.}, 55:213, 1976.

\bibitem{CEMS96}
A.~Calini, N.~M. Ercolani, D.~W. McLaughlin, and C.~M. Schober.
\newblock Melnikov {A}nalysis of {N}umerically {I}nduced {C}haos in the
  {N}onlinear {S}chr{\"{o}}dinger {E}quation.
\newblock {\em Phys. D}, 89, no.3-4:227--260, 1996.

\bibitem{Die60}
J.~Dieudonne.
\newblock {\em Foundations of {M}odern {A}nalysis}.
\newblock Academic Press, 1960.

\bibitem{Fen79}
N.~Fenichel.
\newblock Geometric {S}ingular {P}erturbation {T}heory for {O}rdinary
  {D}ifferential {E}quations.
\newblock {\em J Diff Eqns}, 31:53--98, 1979.

\bibitem{Gra74}
S.~M. Graff.
\newblock On the {C}onservation of {H}yperbolic {I}nvariant {T}ori for
  {H}amiltonian {S}ystems.
\newblock {\em J. Diff. Eqs.}, 15:1--69, 1974.

\bibitem{Li92}
Y.~Li.
\newblock Backlund {T}ransformations and {H}omoclinic {S}tructures for the
  {N}{L}{S} {E}quation.
\newblock {\em Phys. Letters A}, 163:181--187, 1992.

\bibitem{Li99b}
Y.~Li.
\newblock Homoclinic {T}ubes in {N}onlinear {S}chr{\"{o}}dinger {E}quation
  {U}nder {H}amiltonian {P}erturbations.
\newblock {\em Progress of Theoretical Physics}, 101, No. 3(4):559--577, 1999.

\bibitem{Li01d}
Y.~Li.
\newblock On 2{D} {E}uler {E}quations: {P}art {I}{I}. {L}ax {P}airs and
  {H}omoclinic {S}trucutres.
\newblock {\em Submitted}, 2001.

\bibitem{LM94}
Y.~Li and D.~W. McLaughlin.
\newblock Morse and {M}elnikov {F}unctions for {N}{L}{S} {P}des.
\newblock {\em Comm. Math. Phys.}, 162:175--214, 1994.

\bibitem{LM97}
Y.~Li and D.~W. McLaughlin.
\newblock Homoclinic {O}rbits and {C}haos in {P}erturbed {D}iscrete {N}{L}{S}
  {S}ystem. {P}art {I} {H}omoclinic {O}rbits.
\newblock {\em Journal of Nonlinear Sciences}, 7:211--269, 1997.

\bibitem{LW97a}
Y.~Li and S.~Wiggins.
\newblock Homoclinic {O}rbits and {C}haos in {P}erturbed {D}iscrete {N}{L}{S}
  {S}ystem. {P}art {II} {S}ymbolic {D}ynamics.
\newblock {\em Journal of Nonlinear Sciences}, 7:315--370, 1997.

\bibitem{LW97b}
Y.~Li and S.~Wiggins.
\newblock {\em Invariant {M}anifolds and {F}ibrations for {P}erturbed
  {N}onlinear {S}chr{\"{o}}dinger {E}quations}, volume 128.
\newblock Springer-Verlag, Applied Mathematical Sciences, 1997.

\bibitem{MO95}
D.W. McLaughlin and E.~A. Overman.
\newblock Whiskered {T}ori for {I}ntegrable {P}des and {C}haotic {B}ehavior in
  {N}ear {I}ntegrable {P}des.
\newblock {\em Surveys in Appl. Math. 1}, 1993.

\bibitem{Pal88}
K.~J. Palmer.
\newblock Exponential {D}ichotomies, the {S}hadowing {L}emma and {T}ransversal
  {H}omoclinic {P}oints.
\newblock {\em Dynamics Reported}, 1:265--306, 1988.

\bibitem{Sil67}
L.~P. Silnikov.
\newblock The {E}xistence of a {C}ountable {S}et of {P}eriodic {M}otions in the
  {N}eighborhood of a {H}omoclinic {C}urve.
\newblock {\em Soviet Math. Dokl.}, 8:102--106, 1967.

\bibitem{Sil68}
L.~P. Silnikov.
\newblock Structure of the {N}eighborhood of a {H}omoclinic {T}ube of an
  {I}nvariant {T}orus.
\newblock {\em Soviet Math. Dokl.}, 9, No.3:624--628, 1968.

\end{thebibliography}

\end{document}